\documentclass{article}
\usepackage{amsmath}
\usepackage{multicol,color}
\usepackage{float}
\usepackage{soul}
\usepackage{graphicx}
\usepackage{amssymb}
\usepackage{theorem}
\usepackage{fancyhdr}
\usepackage{tikz}
\usepackage{enumerate}
\usepackage[margin=1.1in]{geometry}
\def\disp{\displaystyle}
\def\ve{\varepsilon}

\def\lm{\lambda}
\def\O{\Omega}
\def\Tilde{\widetilde}

\def\tx{\widetilde{x}}
\def\tu{\widetilde{u}}
\def\ta{\widetilde{a}}
\def\tz{\widetilde{z}}

\def\oa{\bar a}

\def\ox{\bar{x}}
\def\oy{\bar{y}}
\def\oz{\bar{z}}

\def\ou{\bar{u}}

\def\cone{\hbox{}}

\def\gph{\hbox{}}

\def\gg{\gamma}
\def\dn{\downarrow}

\def\tto{\rightrightarrows}

\def\Limsup{\mathop{{\rm Lim}\,{\rm sup}}}

\def\Hat{\widehat}
\def\Tilde{\widetilde}
\def\Bar{\overline}
\def\ra{\rangle}
\def\la{\langle}
\def\ve{\varepsilon}

\def\h{\hfill\Box}
\def\R{\mathbb{R}}
\def\N{\mathbb{N}}

\def\gph{\mbox{\rm gph}\,}
\def\epi{\mbox{\rm epi}\,}

\def\dom{\mbox{\rm dom}\,}

\def\dist{\mbox{\rm dist}}

\def\cone{\mbox{\rm cone}\,}

\def\var{\mbox{\rm var}\,}
\def\dn{\downarrow}
\def\O{\Omega}

\def\emp{\emptyset}

\def\oR{\Bar{\R}}
\def\lm{\lambda}

\def\gg{\gamma}

\def\al{\alpha}

\def\N{I\!\!N}
\def\th{\theta}

\newtheorem{theorem}{Theorem}[section]

\newtheorem{proposition}[theorem]{Proposition}

\theoremstyle{plain}{\theorembodyfont{\rmfamily}
}
\theoremstyle{plain}{\theorembodyfont{\rmfamily}
}
\theoremstyle{plain}{\theorembodyfont{\rmfamily}
}
\theoremstyle{plain}{\theorembodyfont{\rmfamily}
\newtheorem{example}[theorem]{Example}}

\theoremstyle{plain}{\theorembodyfont{\rmfamily}
}

\begin{document}
\begin{center}
{\bf OPTIMAL CONTROL OF A PERTURBED SWEEPING PROCESS\\VIA DISCRETE APPROXIMATIONS}\footnote{This research was partly supported by the National Science Foundation under grants DMS-1007132 and DMS-1512846 and by the Air Force Office of Scientific Research grant \#15RT0462.}\\[3ex]
TAN H. CAO\footnote{Department of Mathematics, Wayne State University, Detroit, Michigan 48202, USA (tan.cao@wayne.edu).} and BORIS S. MORDUKHOVICH\footnote{Department of Mathematics, Wayne State University, Detroit, Michigan 48202, USA (boris@math.wayne.edu).}
\end{center}
\vspace*{0.05in}
\small{\sc Abstract.} The paper addresses an optimal control problem for a perturbed sweeping process of the rate-independent hysteresis type described by a controlled ``play-and stop" operator with separately controlled perturbations. This problem can be reduced to dynamic optimization of a state-constrained unbounded differential inclusion with highly irregular data that cannot be treated by means of known results in optimal control theory for differential inclusions. We develop the method of discrete approximations, which allows us to adequately replace the original optimal control problem by a sequence of well-posed finite-dimensional optimization problems whose optimal solutions strongly converge to that of the controlled perturbed sweeping process. To solve the discretized control systems, we derive effective necessary optimality conditions by using second-order generalized differential tools of variational analysis that explicitly calculated in terms of the given problem data.\\
{\em Key words and phrases.} Controlled sweeping process, Play-and-stop operator, Hysteresis, Discrete approximations, Variational analysis, Generalized differentiation, Optimality conditions.\\
{\em 2010 Mathematics Subject Classification.} Primary: 49M25, 47J40; Secondary: 90C30, 49J53.\vspace*{-0.1in}

\normalsize
\section{Introduction and Problem Formulation}
\setcounter{equation}{0}\vspace*{-0.1in}

In this paper we deal with a version of the {\em sweeping process} introduced by Jean-Jacques Moreau in the 1970s (see his comprehensive for that time lecture notes \cite{mor_frict} with the references to the original publications) in the following form of the dissipative differential inclusion:
\begin{equation}\label{e:1}
-\dot{x}(t)\in N\big(x(t);C(t)\big)\;\mbox{ a.e. }\;t\in[0,T],\quad x(0):=x_0\in C(0)\subset H,
\end{equation}
where $C(t)$ is a continuously moving convex set, and where the normal cone operator to a convex subset $C\subset H$ of a Hilbert space is given by
\begin{eqnarray}\label{nor}
N(x;C):=\big\{v\in H\big|\;\la v,y-x\ra\le 0,\;y\in C\big\}\;\mbox{ if }\;x\in C\;\mbox{ and }\;N(x;C):=\emp\;\mbox{ if }\;x\notin C.
\end{eqnarray}
The latter construction allows us to equivalently describe \eqref{e:1} as an {\em evolution variational inequality} \cite{BrS,Kr}, or as a {\em differential variational inequality} in the terminology of \cite{ps,St}. The original motivations for the introduction and study of the sweeping process come from applications to mechanical systems mostly related to friction and elastoplasticity, while further developments apply also to various problems of hysteresis, ferromagnetism, electric circuits, phase transitions, economics, etc.; see, e.g., \cite{aht,BrS,KP,Kr,MM,smb,St} and the extensive bibliographies therein.\vspace*{-0.05in}

It has been realized in the sweeping process theory \cite{CT,KMM} that the Cauchy problem \eqref{e:1} admits a {\em unique} solution under mild assumptions on the given moving set $C(t)$. This excludes considering any optimization problem for \eqref{e:1} and thus strictly distinguishes the sweeping process from the conventional optimal control theory for differential inclusions $\dot x(t)\in F(x(t))$. The latter theory supposes the existence of numerous solutions to the Cauchy problem and then studies minimization of some cost functionals over them; see, e.g., \cite{m-book2,v} with further references and discussions. \vspace*{-0.05in}

There are three approaches in the literature to introduce control actions in the sweeping process frameworks and then to conduct optimization with respect to these controls and the corresponding sweeping trajectories. The {\em first approach} considers controls in additive {\em perturbations} on the right-hand side of \eqref{e:1} without changing the moving set $C(t)$. The results obtained in this direction mostly concern existence theorems and relaxation procedures while not optimality conditions; see \cite{et} and the recent papers \cite{aht,cmf} with the references therein. We also place into this category the main results of \cite{ri1,ri2} that establish the existence of optimal controls for rate-independent evolutions via some direct methods involving finite-difference approximations. The {\em second approach} developed in \cite{bk} for $H=\R^n$ and then partly extended in \cite{ao} introduces controls in an ordinary differential equation associated with the sweeping process over a {\em given} set $C(t)\subset\R^n$. The obtained results provide necessary optimality conditions for the continuous-time problem in \cite{bk} and for the approximating finite-difference systems in \cite{ao}.\vspace*{-0.05in}

The {\em third approach} to optimal control of the sweeping process, which is implemented in this paper, employs a control parametrization {\em directly in the sweeping set} $C(t)$ making it dependent on control actions. It has been initiated in \cite{chhm1} for the case of a controlled hyperplane  $C(t)$ in $\R^n$ and then strongly developed in the very recent paper \cite{chhm2} for the case of moving controlled polyhedra
\begin{eqnarray}\label{sw-poly}
C(t):=\big\{x\in\R^n\big|\;\la u_i(t),x\ra\le b_i(t),\;i=1,\ldots,m\big\},
\end{eqnarray}
where the control functions $u_i(t)$ and $b_i(t)$ are absolutely continuous on the fixed time interval $[0,T]$. Necessary optimality conditions for such optimal control problems of a new type are derived in \cite{chhm1,chhm2} by using discrete approximations and appropriate generalized differential tools of variational analysis in the lines of \cite{m95,m-book1,m-book2}. Note that the discrete approximation approach to optimization of bounded differential inclusions $\dot x\in F(x)$ developed in \cite{m95,m-book2} essentially relies on the Lipschitz continuity of $F$, which fails in the case of \eqref{e:1} with controlled sweeping sets as in \eqref{sw-poly} and thus requires serious modifications.

This paper concerns a new class of optimal control problems for the {\em perturbed sweeping process}
\begin{equation}\label{e:5}
 -\dot x(t)\in N\big(x(t);C(t)\big)+f\big(x(t),a(t)\big)\;\mbox{ a.e. }\;t\in[0,T],\quad x(0):=x_{0}\in C(0),
\end{equation}
with one part of controls $a\colon[0,T]\to\R^d$ acting in the perturbation mapping $f\colon\R^n\times\R^d\to\R^n$ and the other part of controls $u\colon[0,T]\to\R^n$ acting in the moving set
\begin{equation}\label{e:6}
C(t):=C+u(t)\;\mbox{ with }\;C:=\big\{x\in\R^n\big|\;\la x^*_i,x\ra\le 0\;\mbox{ for all }\;i=1,\ldots,m\big\},
\end{equation}
where $x^*_i$ are fixed vectors from $\R^n$, and where the final time $T>0$ is also fixed. The minimizing cost functional is given in the generalized Bolza form
\begin{equation}\label{e:4}
\mbox{minimize }\;J[x,u,a]:=\varphi\big(x(T)\big)+\int_0^{T}\ell\big(t,x(t),u(t),a(t),\dot x(t),\dot u(t),\dot a(t)\big)\,dt
\end{equation}
with the proper terminal extended-real-valued cost function $\varphi\colon\R^n\to\oR:=(-\infty,\infty]$ and the running cost function $\ell:[0,T]\times\R^{4n+2d}\to\oR$.
Fixed $r>0$, we impose the additional constraint on the $u$-controls:
\begin{equation}\label{e:7}
\|u(t)\|=r\;\mbox{ for all }\;t\in[0,T]
\end{equation}
required by applications. Note that the controlled moving set in \eqref{e:6} can be written in the polyhedral form $C(t)=\{x\in\R^n|\;\la x^*_i,x\ra\le-u_i(t),\;i=1,\ldots,m\}$ for $u(t)=(u_1(t),\ldots,u_m(t)$, which shows that \eqref{e:6} reduces to \eqref{sw-poly} with controls acting only on the right-hand sides of the polyhedral inequalities. However, this reduction does not allow us to apply the results of \cite{chhm2} to our problem, even in the absence of perturbations, since there are no constraints on the controls $b_i(t):=-u_i(t)$ in \cite{chhm2}, while the imposed control constraint \eqref{e:7} cannot be ignored in our setting.\vspace*{-0.05in}

Besides the constraint issue, a major ingredient that distinguishes the novel framework in \eqref{e:5}, \eqref{e:6} from the one in \eqref{e:1}, \eqref{sw-poly} is the {\em presence of controls} in {\em perturbations} together with those in the moving {\em sweeping set}. This makes the new model challenging from the viewpoint of variational analysis and important for various applications. The primary application we have in mind is the {\em crowd motion model} (see, e.g., \cite{mv}), which corresponds to \eqref{e:5} with controls only in perturbations and whose simplified optimal control version is solved in our adjacent paper \cite{cm1} based on the obtained optimality conditions.\vspace*{-0.05in}

From a different viewpoint, the given description \eqref{e:6} of the moving set in \eqref{e:5} for each fixed function $u(\cdot)$ relates to the so-called {\em play-and-stop operator}; see \cite{KP,Kr,smb} and the references therein. As discussed in \cite{smb}, such operators constitute basic elements of the mathematical theory of {\em rate independent hysteresis processes} including, in particular, the celebrated Preisach model in ferromagnetism.\vspace*{-0.05in}

Yet another essential difference between the frameworks of \cite{chhm2} and of the current paper, even in the absence of controlled perturbations, is the choice of classes of {\em feasible controls}. The most natural class in the setting of \cite{chhm2} is the collection of controls $u(t)$ absolutely continuous on $[0,T]$, which generate absolutely continuous trajectories $x(t)$ of \eqref{sw-poly} under appropriate qualification conditions; cf.\ \cite{chhm2} for more details. In the setting \eqref{e:5}, \eqref{e:6} of this paper by feasible controls $u(\cdot)$ and $a(\cdot)$ we understand functions that belong to the Sobolev spaces $W^{1,2}([0,T];\R^n)$ and $W^{1,2}([0,T];\R^n)$, respectively. It follows from the powerful well-posedness result of \cite{et} that such a control pair generates a unique trajectory $x(\cdot)\in W^{1,2}([0,T];\R^n)$ of the sweeping inclusion \eqref{e:5}. Having this in mind, we formulate the {\em sweeping optimal control problem} $(P)$ as follows: minimize \eqref{e:4} over $W^{1,2}$-controls $(u(\cdot),a(\cdot))$ on $[0,T]$ and the corresponding $W^{1,2}$-trajectories $x(\cdot)$ of \eqref{e:5} with $C(t)$ from \eqref{e:6} subject to the control equality constraint \eqref{e:7}. Observe that, besides \eqref{e:7}, there are implicit {\em mixed} (i.e., state-control) inequality constraints in $(P)$ given by
\begin{equation}\label{mixed}
\big\la x^*_i,x(t)-u(t)\big\ra\le 0\;\mbox{ for all }\;t\in[0,T]\;\mbox{ and }\;i=1,\ldots,m,
\end{equation}
which follow from \eqref{e:5} and \eqref{e:6} due to the second part of the normal cone definition \eqref{nor}.\vspace*{-0.05in}

The {\em main goal} of this paper is to study the formulated optimal control problem $(P)$ and its slight parametric modification $(P^\tau)$ defined below by using the {\em method of discrete approximations} in the vein of \cite{m95,m-book2} and its significant elaboration for the case of unperturbed non-Lipschitzian differential inclusions developed in \cite{chhm2}. The presence of controlled perturbations in \eqref{e:5} together with the control and mixed constraints in \eqref{e:7},\eqref{mixed} essentially complicates the discrete approximation procedure. We aim to construct well-posed discrete approximations in such a way that every {\em feasible} (resp.\ {\em locally optimal}) solution to $(P^\tau)$ with $\tau\ge 0$ and $(P^0)=(P)$, can be {\em strongly approximated} in $W^{1,2}[0,T]$ by feasible (resp.\ optimal) solutions to finite-difference control systems. Employing then appropriate first-order and second-order generalized differential constructions of variational analysis and explicitly calculating them via the problem data allow us to obtain effective {\em necessary optimality conditions} for discrete optimal solutions, which can be treated as {\em suboptimality} (almost optimality) conditions for the original sweeping control problem. Deriving exact optimality conditions for the continuous-time sweeping process in $(P^\tau)$ is a subject of \cite{cm1}. Note that the finite-difference problems constructed in this paper provide more precise approximations of feasible and optimal solutions of $(P^\tau)$ in comparison with the corresponding settings of \cite{chhm2} and lead us in this way to new conditions for optimality in both discrete and continuous frameworks.\vspace*{-0.05in}

The rest of the paper is organized as follows. Section~2 lists basic assumptions and preliminaries from the sweeping process theory used below. In Section~3 we justify the existence of optimal solutions to $(P^\tau)$ as $\tau\ge 0$ and discuss its relaxation stability. This section also contains the definition of intermediate local minimizers and their relaxed modification studied in the paper.\vspace*{-0.05in}

In Section~4 we develop a constructive finite-difference procedure to strongly in $W^{1,2}[0,T]$ approximate any feasible control $(u(t),a(t))$ and the corresponding trajectory $x(t)$ of \eqref{e:5} by feasible solutions to discrete approximation systems that are piecewise linearly extended to $[0,T]$. The next Section~5 establishes the strong $W^{1,2}$-approximation of the local optimal solution to $(P^\tau)$ by extended optimal solutions to its discrete counterparts constructed therein. This makes a bridge between optimization of the continuous-time and discrete-time sweeping control systems while justifying in this way an effective usage of discrete approximations to solve the original sweeping control problem.\vspace*{-0.05in}

After reviewing in Section~6 the generalized differential tools of variational analysis used in the paper and their calculations via the given data of $(P^\tau)$, we derive in Section~7 necessary optimality conditions for the constructed discrete approximations. Section~8 presents a numerical example of using these conditions to find optimal solutions. Finally, in Section~9 we discuss the possibility to employ the obtained results in the limiting procedure to derive nondegenerate necessary optimality conditions for a given local minimizer in $(P^\tau)$ with  their subsequent applications to the crowd motion model.\vspace*{-0.05in}

The notation of this paper is standard in variational analysis and optimal control; see, e.g., \cite{m-book1,rw,v}. Recall that the symbol $B(x,\ve)$ denotes the closed ball of the space in question centered at $x$ with radius $\ve>0$ while $\N$ signifies the collection of all natural numbers $\{1,2,\ldots\}$.\vspace*{-0.15in}

\section{Standing Assumptions and Preliminaries}
\setcounter{equation}{0}\vspace*{-0.1in}

Throughout the paper we impose the following standing assumptions on the initial data of the optimal control problem $(P)$ in \eqref{e:5}--\eqref{mixed}:

{\bf (H1)} The mapping $f\colon\R^n\times\R^d\to\R^n$ is continuous on $\R^n\times\R^d$ and locally Lipschitz continuous in the first argument, i.e.,
for every $\varepsilon>0$ there is a constant $K>0$ such that
\begin{equation}\label{f-lip}
\left\|f(x,a)-f(y,a)\right\|\le K\left\|x-y\right\|\;\mbox{ whenever }\;(x,y)\in B(0,\varepsilon)\times B(0,\varepsilon),\quad a\in\R^d.
\end{equation}
Furthermore, there is a constant $M>0$ ensuring the growth condition
\begin{equation}\label{gro}
\left\|f(x,a)\right\|\le M\big(1+\|x\|\big)\;\mbox{ for any }\;x\in\bigcup_{t\in[0,T]}C(t),\quad a\in\R^d.
\end{equation}\vspace*{-0.1in}

{\bf (H2)} The terminal cost function $\varphi:\R^n\to\Bar{\R}$ and the running cost function $\ell:[0,T]\times\R^{4n+2d}\to\Bar{\R}$ in \eqref{e:4} are lower semicontinuous (l.s.c.) while $\ell$ is bounded from below on bounded sets.

Now we are ready to formulate the powerful well-posedness result for the sweeping process under consideration that reduces to \cite[Theorem~1]{et}.\vspace*{-0.1in}

\begin{proposition}{\bf (well-posedness of the controlled sweeping process).}\label{Th:1} Under the assumptions in {\rm(H1)}, let $u(\cdot)\in W^{1,2}([0,T];\R^n)$ and $a(\cdot)\in W^{1,2}([0,T];\R^d)$, and let $M>0$ be taken from \eqref{gro} Then the perturbed sweeping inclusion \eqref{e:5} with $C(t)$ from \eqref{e:6} admits the unique solution $x(\cdot)\in W^{1,2}([0,T];\R^n)$ generated by $(u(\cdot),a(\cdot))$ and satisfying the estimates
\begin{eqnarray*}
\|x(t)\|\le l:=\|x_0\|+e^{2MT}\Big(2MT(1+\|x_0\|)+\int_0^T\|\dot u(s)\|ds\Big)\;\mbox{ for all }\;t\in[0,T],
\end{eqnarray*}
\begin{equation*}
\|\dot{x}(t)\|\le 2(1+l)M+\|\dot u(t)\|\;\mbox{ a.e. }\;t\in[0,T].
\end{equation*}
\end{proposition}\vspace*{-0.1in}
{\bf Proof.} To deduce this result from \cite[Theorem~1]{et}, with taking into account the solution estimates therein, it remains to verify that $C(t)$ in \eqref{e:6} generated by the chosen $W^{1,2}$-control $u(\cdot)$ varies in an {\em absolutely continuous way} \cite{et}, i.e., there is an absolutely continuous function $v\colon[0,T]\to\R$ such that
\begin{equation}\label{abs}
\left|\dist\big(y;C(t)\big)-\dist\big(y;C(s)\big)\right|\le\left|v(t)-v(s)\right|\;\mbox{ for all }\;t,s\in[0,T]
\end{equation}
with $\dist(x;\O)$ standing for the distance from $x\in\R^n$ to the closed set $\O\subset\R^n$ and with the function
$$
v(t):=\int^t_0\|\dot u(s)\|ds,\quad 0\le t\le T,
$$
in our case. To verify \eqref{abs}, pick any $y\in\R^n$ and $c\in C$ and then easily get the estimates
$$
\dist\big(y;C(t)\big)=\dist\big(y;u(t)+C\big)\le\left\|y-u(t)-c\right\|\le\left\|y-u(s)-c\right\|+\left\|u(t)-u(s)\right\|,
$$
which imply in turn by the definition of the distance function that
$$
\dist\big(y;C(t)\big)\le\inf_{c\in C}\left\|y-u(s)-c\right\|+\left\|u(t)-u(s)\right\|=\dist\big(y;C(s)\big)+\left\|u(t)-u(s)\right\|.
$$
Using this and then changing the positions of $t$ and $s$ give us the resulting inequality
$$
\left|\dist\big(y;C(t)\big)-\dist\big(y;C(s)\big)\right|\le\left\|u(t)-u(s)\right\|\;\mbox{ for all }\;t,s\in[0,T].
$$
This finally yields \eqref{abs} by observing that
$$
\left|d\big(y;C(t)\big)-d\big(y;C(s)\big)\right|\le\left\|u(t)-u(s)\right\|=\left\|\int^t_s\dot{u}(\th)d\th\right\|\le\int^t_s\|\dot{u}(\th)\|d\th=
\int^t_s\dot{v}(\th)d\th=|v(t)-v(s)|
$$
and thus completes the proof of the proposition.$\h$\vspace*{-0.15in}

\section{Discrete Approximations of Feasible Solutions}
\setcounter{equation}{0}\vspace*{-0.1in}

In this section we construct a sequence of discrete approximations of the sweeping differential inclusion in \eqref{e:5}, \eqref{e:6} with the constraints in \eqref{e:7} and \eqref{mixed}, but without appealing to the minimizing functional \eqref{e:4}. The main result of this section is justifying the strong $W^{1,2}$-approximation of {\em any} feasible control and the corresponding sweeping trajectory by their finite-difference counterparts, which are piecewise linearly extended to the continuous-time interval $[0,T]$.\vspace*{-0.05in}

First we reduce \eqref{e:5} to a more conventional form of differential inclusions. Introduce the new variable $z:=(x,u,a)\in\R^n\times\R^n\times\R^d$ and define the set-valued mapping $F:\R^n\times\R^n\times\R^d\tto\R^n$ by
\begin{equation}\label{e:17}
F(z)=F(x,u,a):=N(x-u;C)+f(x,a).
\end{equation}
Consider the collection of {\em active constraint indices} of polyhedron \eqref{e:6} at $\ox\in C$ given by
\begin{equation}\label{e:19}
I(\bar{x}):=\big\{i\in\{1,\ldots,m\}\big|\;\big\la x^*_i,\bar{x}\big\ra=0\big\},
\end{equation}
it is not difficult to observe (see, e.g., \cite[Proposition~3.1]{hmn}) the explicit representation
\begin{equation}\label{e:20}
F(z)=\Big\{\sum_{i\in I(x-u)}\lambda_{i}x^*_i\Big|\;\lambda_i\ge 0\Big\}+f(x,a)
\end{equation}
of \eqref{e:17} via the active index set \eqref{e:19} at $x-u\in C$. Then we can rewrite \eqref{e:5} in the following equivalent form with respect to the variable $z=(x,u,a)$:
\begin{equation}\label{e:21}
-\dot{z}(t)\in F\big(z(t)\big)\times\R^n\times\R^d\;\mbox{ a.e. }\;t\in[0,T]
\end{equation}
with the initial condition $z(0)=(x_{0},u(0),a(0))$ satisfying $x_{0}-u(0)\in C$, i.e., such that $\la x^*_i,x_0-u(0)\ra\le 0$ for all $i=1,\ldots,m$. Proposition~\ref{Th:1} allows us to have solutions of the differential inclusion \eqref{e:21} in the class of $W^{1,2}$-functions $z(t)=(x(t),u(t),a(t))$ on $[0,T]$.\vspace*{-0.05in}

Note that, although the resulting system \eqref{e:21} is written in the conventional form of the theory of differential inclusions, it does not satisfy usual assumptions therein. Indeed, the right-hand side of \eqref{e:21} is intrinsically {\em unbounded} in all its components, including the first (perturbed normal cone) one in which is {\em highly non-Lipschitzian}. Furthermore, the constrained system under consideration contains the intrinsic inequality {\em state constraints} \eqref{mixed} together with the equality one \eqref{e:7} on the whole time interval $[0,T]$.\vspace*{-0.05in}

Having in mind further applications including those developed in \cite{cm1}, it makes sense to consider a {\em parametric version} of the equality constraint in \eqref{e:7} with a small parameter $\tau\ge 0$ while replacing \eqref{e:7} by
\begin{eqnarray}\label{e:8}
\left\{\begin{array}{ll}
\|u(t)\|=r\;\mbox{ for all }\;t\in[\tau,T-\tau],\\
r-\tau\le\|u(t)\|\le r+\tau\;\mbox{ for all }\;t\in[0,\tau)\cup(T-\tau,T],
\end{array}\right.
\end{eqnarray}
which reduces to \eqref{e:7} when $\tau=0$. Fix any $\tau\in[0,\min\{r,T\}]$, $k\in\N$ and denote by $j_\tau(k):=\left[k\tau/T\right]$ the smallest index $j$ such that $t^k_j\ge\tau$ and by $j^\tau(k):=[k(T-\tau)/T]-1$ the largest $j$ with $t^k_j\le T-\tau$.\vspace*{-0.05in}

The next theorem on the {\em strong discrete approximation} of feasible sweeping solutions is a counterpart of \cite[Theorem~3.1]{chhm2} for the perturbed sweeping process in \eqref{e:5}, \eqref{e:6} constrained by \eqref{mixed}, \eqref{e:8} with additional {\em quantitative estimates} expressed via the system data. The reader can see that both the formulation and proof in the new setting are significantly more involved in comparison with \cite{chhm2}. Observe also the {\em novel approximation conclusion} \eqref{x0}, which holds also in the setting of \cite{chhm2} while being missed therein. This conclusion will allow us to construct a more precise discrete approximation of a local minimizer in Theorem~\ref{Th:5}, which is crucial to derive a new transversality condition for the original continuous-time sweeping control problem $(P^\tau)$ in \cite{cm1}.\vspace*{-0.1in}

\begin{theorem}{\bf($W^{1,2}$-strong discrete approximation of feasible sweeping solutions).}\label{Th:3} Under the validity of {\rm(H1)}, let the triple $\oz(\cdot)=(\ox(\cdot),\ou(\cdot),\oa(\cdot))$ be a given feasible solution to the constrained sweeping system from \eqref{e:5}, \eqref{e:6}, and \eqref{e:8} with a fixed parameter $\tau\in[0,\min\{r,T\}]$, and let the constant $K$ be taken from \eqref{f-lip}. Define the discrete partitions of $[0,T]$ by
\begin{equation}\label{mesh}
\Delta_{k}:=\big\{0=t^{k}_{0}<t^{k}_{1}<\ldots<t^{k}_{k}\big\}\;\mbox{ with }\;h_{k}:=t^{k}_{j+1}-t^{k}_{j}\dn 0\;\mbox{ as }\;k\to\infty
\end{equation}
and suppose that $\oz(\cdot)$ has the following properties at the mesh points $($while observing that all these properties hold if $\oz(\cdot)\in W^{2,\infty}[0,T])$: it satisfies \eqref{e:21} at $t^k_j$ as $j=0,\ldots,k-1$ for all $k\in\N$ $($with the right-side derivative at $t_0=0)$, we have
\begin{equation}\label{e:25}
\sum^{k-1}_{j=0}(t^k_{j+1}-t^k_j)\left\|\dfrac{\ox(t^k_{j+1})-\ox(t^k_j)}{t^k_{j+1}-t^k_j}-\dot{\ox}(t^k_j)\right\|^2\to 0\;\mbox{ as }\;k\to\infty,
\end{equation}
and there is a constant $\mu>0$ independent of $k$ such that
\begin{equation}\label{e:26}
\begin{array}{ll}
\disp\sum^{k-1}_{j=0}\left\|\dfrac{\ox(t^k_{j+1})-\ox(t^k_j)}{t^k_{j+1}-t^k_j}-\dot{\ox}(t^k_j)\right\|\le\mu,\quad\left\|\dfrac{\ou(t^k_1)-\ou(t^k_0)}{t^k_1-t^k_0}
\right\|\le\mu,\\\
\disp\sum^{k-2}_{j=0}\left\|\dfrac{\ou(t^k_{j+2})-\ou(t^k_{j+1})}{t^k_{j+2}-t^k_{j+1}}
\disp-\dfrac{\ou(t^k_{j+1})-\ou(t^k_j)}{t^k_{j+1}-t^k_j}\right\|\le\mu.
\end{array}
\end{equation}
Then there exist a sequence of piecewise linear functions $z^{k}(t):=(x^{k}(t),u^{k}(t),a^{k}(t))$ on $[0,T]$ and a sequence of $\ve_k\le 2h_k\mu e^K\dn 0$ as $k\to\infty$ for which $(x^k(0),u^k(0),a^k(0))=(x_0,\ou(0),\oa(0))$,
\begin{equation}\label{x0}
\frac{x^k(t^k_1)-x^k(t^k_0)}{h_k}\to\dot\ox(0)\;\mbox{ as }\;k\to\infty,
\end{equation}
\begin{equation}\label{e:28}
\left\{\begin{array}{lcl}
\|u^{k}(t^{k}_{j})\|=r\;\mbox{ if }\;j=j_\tau(k),\ldots,j^\tau(k),\\
r-\tau-\varepsilon_k\le\|u^k(t^k_j)\|\le r+\tau+\varepsilon_k\;\mbox{ if }\;j=0,\ldots,j_\tau(k)-1\;\mbox{ and }\;j\ge j^\tau(k)+1,
\end{array}\right.
\end{equation}
\begin{equation}\label{e:29}
x^{k}(t)=x^{k}(t_{j})-(t-t_{j})v^{k}_{j},\;x^k(0)=x_{0},\;t^{k}_{j}\le t\le t^{k}_{j+1}\;\mbox{ with }\;v^{k}_{j}\in F\big(z^{k}(t^{k}_{j})\big),\quad j=0,\ldots,k-1,
\end{equation}
and the functions $z^k(\cdot)$ converge to $\oz(\cdot)$ in the norm topology of $W^{1,2}[0,T]$, i.e.,
\begin{equation}\label{e:30}
z^{k}(t)\to\bar{z}(t)\;\mbox{ uniformly on }\;[0,T]\;\mbox{ and }\;\int_{0}^{T}\|\dot{z}^{k}(t)-\dot{\bar{z}}(t)\|^{2}dt\to 0\;\mbox{ as }\;k\to \infty.
\end{equation}
Furthermore, for every $k\in\N$ we have the estimates
\begin{equation}\label{e:31}
\var\big(\dot u^k;[0,T]\big)\le\Tilde{\mu}\;\mbox{ and }\;\left\|\dfrac{u^k(t^k_1)-u^k(t^k_0)}{h_k}\right\|\le\Tilde{\mu}\;\mbox{ with }\;\Tilde{\mu}:=\max\big\{3\mu(1+4K)e^K,4\mu(e^K+1)\big\},
\end{equation}
where the symbol $``\var"$ stands for the total variation of the function in question.
\end{theorem}\vspace*{-0.1in}
{\bf Proof.} Let $y^k(\cdot):=(y^k_1(\cdot),y^k_2(\cdot),y^k_3(\cdot))$ be piecewise linear on $[0,T]$ and such that
\begin{equation}\label{yk}
\left(y^k_1(t^k_j),y^k_2(t^k_j),y^k_3(t^k_j)\right):=\left(\ox(t^k_j),\ou(t^k_j),\oa(t^k_j)\right),\quad j=0,\ldots,k.
\end{equation}
Define next $w^k(t)=(w^k_1(t),w^k_2(t),w^k_3(t)):=\dot{y}^k(t)$ as piecewise constant and right continuous function on $[0,T]$ via the derivatives at non-mesh points and deduce from \eqref{e:28} that $\var(w^k_2;[0,T])\le\mu$ for every $k\in\N$. It follows from the definition of $w^k(\cdot)$ that
\begin{equation}\label{w0}
w^k_1(0)=\frac{\ox(t^k_1)-\ox(t^k_0)}{h_k}\to\dot\ox(0)\;\mbox{ as }\;k\to\infty
\end{equation}
due to the existence of the right derivative of $\dot\ox(0)$ by the imposed assumption on the validity of \eqref{e:21} at the mesh points. Furthermore, we get from \eqref{mixed} that
$$
\left\langle x^*_i,y^k_1(t^k_j)-y^k_2(t^k_j)\right\rangle=\left\langle x^*_i,\ox(t^k_j)-\ou(t^k_j)\right\rangle\le 0
$$
on the mesh $\Delta_k$ for all $j=1,\ldots,k-1$ and $i=1,\ldots,m$. The constructions made ensure that
$$
y^k(\cdot)\to\oz(\cdot)\;\mbox{ uniformly on }\;[0,T]\;\mbox{ and }\;w^k(\cdot)\to\dot\oz(\cdot)\;\mbox{ in norm of }\;L^2([0,T];\R^{2n+d}).
$$
Denote $a^k(t):=y^k_3(t)$ for all $t\in[0,T]$, fix $k\in\N$, and use for simplicity the notation $t_j:=t^k_j$ as $j=1,\ldots,k$. To construct the claimed trajectories $x^k(t)$ of \eqref{e:29}, we proceed by induction and suppose that the value of $x^k(t_j)$ is known. Define now the vectors
$$
u^k(t_j):=x^k(t_j)-y^k_1(t_j)+y^k_2(t_j)=x^k(t_j)-\ox(t_j)+\ou(t_j),\quad j=0,\ldots,k,
$$
and assume without loss of generality that $\left\|u^k(t_j)\right\|=r$ for $j=j_\tau(k),\ldots,j^\tau(k)$, which clearly yields
\begin{equation}\label{e:33}
x^k(t_j)-u^k(t_j)=\ox(t_j)-\ou(t_j)\;\mbox{ for }\;j=0,\ldots,k.
\end{equation}
Since the sets $F(z)$ in \eqref{e:17} are closed and convex, we select the unique projection
\begin{equation}\label{e:34}
v^k_j:=\Pi(-w^k_{1j};F\big(x^k(t_j),u^k(t_j),a^k(t_j))\big),\quad j=0,\ldots,k,
\end{equation}
and deduce from \eqref{w0} that $v^k_0\to\dot\ox(0)$ as $k\to\infty$. Defining next $x^k(t):=x^k(t_j)-(t-t_j)v^k_j$ for all $t\in [t_j,t_{j+1}]$ and $j=0,\ldots,k$ shows that the inclusions in \eqref{e:29} are fulfilled and condition \eqref{x0} holds. Furthermore, we deduce from \eqref{e:20} and \eqref{e:33} that
\begin{equation}\label{e:35}
F\big(x^k(t_j),u^k(t_j),a^k(t_j)\big)=F\big(\ox(t_j),\ou(t_j),\oa(t_j)\big)+f\big(x^k(t_j),\oa(t_j)\big)-f\big(\ox(t_j),\oa(t_j)\big)
\end{equation}
at the mesh points. To verify that the triples $(x^k(t),u^k(t),a^k(t))$, $k\in\N$, constructed above satisfy all the conclusions of the theorem, let us first show that
\begin{equation}\label{e:39}
\varepsilon_k:=\left\|x^k(t_j)-\ox(t_j)\right\|\le\max\big\{h_k\mu(1+h_kK),2h_k\mu e^K\big\}=2h_k\mu e^K\;\mbox{ for all }\;j=0,\ldots,k.
\end{equation}
Indeed, picking any $t\in[t_j,t_{j+1}]$ for $j=0,\ldots,k-1$, we have the representation
\begin{eqnarray*}
x^k(t)-y^k_1(t)=x^k(t_j)-\ox(t_j)+(t-t_j)(-v^k_j-w^k_{1j}),
\end{eqnarray*}
which implies in turn the estimate
\begin{eqnarray*}
\begin{split}
\left\|x^k(t)-y^k_1(t)\right\|\le&\left\|x^k(t_j)-\ox(t_j)\right\|+(t_{j+1}-t_j)\left\|-v^k_j-w^k_{1j}\right\|\\
=&\left\|x^k(t_j)-\ox(t_j)\right\|+(t_{j+1}-t_j)\dist\left(-\dfrac{\ox(t_{j+1})-\ox(t_j)}{t_{j+1}-t_j};F\big(z^k(t_j)\big)\right).
\end{split}
\end{eqnarray*}
It then follows from \eqref{e:35} that
\begin{eqnarray*}
\left\|x^k(t)-y^k_1(t)\right\|\le\left\|x^k(t_j)-\ox(t_j)\right\|+h_k\left\|f\big(x^k(t_j),\oa(t_j)\big)-f\big(\ox(t_j),\oa(t_j)\big)\right\|+h_k\left\|
\dfrac{\ox(t_{j+1})-\ox(t_j)}{t_{j+1}-t_j}-\dot\ox(t_j)\right\|.
\end{eqnarray*}
Using the Lipschitz continuity of $f$ with respect to $x$ imposed in \eqref{f-lip} gives us
\begin{equation}\label{e:36}
\left\|x^k(t)-y^k_1(t)\right\|\le(1+h_kK)\left\|x^k(t_j)-\ox(t_j)\right\|+h_k\left\|\dfrac{\ox(t_{j+1})-\ox(t_j)}{t_{j+1}-t_j}-\dot\ox(t_j)\right\|,
\end{equation}
and thus, by taking the first condition in \eqref{e:26} into account, we arrive at the inequalities
\begin{eqnarray}\label{e:37}
\left\{\begin{array}{ll}
\left\|x^k(t_1)-\ox(t_1)\right\|\le h_k\left\|\dfrac{\ox(t_1)-\ox(t_0)}{t_1-t_0}-\dot\ox(t_0)\right\|,\\
\left\|x^k(t_2)-\ox(t_2)\right\|\le(1+h_kK)\left\|x^k(t_1)-\ox(t_1)\right\|+h_k\left\|\dfrac{\ox(t_2)-\ox(t_1)}{t_2-t_1}-\dot\ox(t_1)\right\|\\
\le h_k\left(\left\|\dfrac{\ox(t_1)-x(t_0)}{t_1-t_0}-\dot x(t_{0})\right\|+\left\|\dfrac{\ox(t_2)-\ox(t_1)}{t_2-t_1}-\dot\ox(t_1)\right\|\right)+h^2_kK\mu.
\end{array}\right.
\end{eqnarray}
Now we proceed by induction to verify that
\begin{equation}\label{e:38}
\left\|x^k(t_j)-\ox(t_j)\right\|\le h_k\sum^{j-1}_{i=0}\left\|\dfrac{\ox(t_{i+1})-\ox(t_i)}{t_{i+1}-t_i}-\dot\ox(t_i)\right\|+
h_k^2K\mu\sum_{i=0}^{j-3}(1+h_kK)^i+(1+h_kK)^{j-1}h_k\mu
\end{equation}
for $j=3,\ldots,k$. Starting with $j=3$, observe from \eqref{e:26}, \eqref{e:36}, and \eqref{e:37} that
\begin{eqnarray*}
\begin{split}
\left\|x^k(t_3)-\ox(t_3)\right\|\le&(1+h_kK)\left\|x^k(t_2)-\ox(t_2)\right\|+h_k\left\|\dfrac{\ox(t_3)-\ox(t_2)}{t_3-t_2}-\dot\ox(t_2)\right\|\\
\le&(1+h_kK)\bigg(h_k\sum^1_{i=0}\left\|\dfrac{\ox(t_{i+1})-\ox(t_i)}{t_{i+1}-t_i}-\dot\ox(t_i)\right\|+h^2_kK\mu\bigg)+h_k\left\|\dfrac{\ox(t_3)-
\ox(t_2)}{t_3-t_2}-\dot\ox(t_2)\right\|\\=&h_k\sum^2_{i=0}\left\|\dfrac{\ox(t_{i+1})-\ox(t_i)}{t_{i+1}-t_i}-\dot\ox(t_i)\right\|+h_k^2K\sum^1_{i=0}
\left\|\dfrac{\ox(t_{i+1})-\ox(t_i)}{t_{i+1}-t_i}-\dot\ox(t_i)\right\|+(1+h_kK)h_k^2K\mu\\
\le&h_k\sum^2_{i=0}\left\|\dfrac{\ox(t_{i+1})-\ox(t_i)}{t_{i+1}-t_i}-\dot\ox(t_i)\right\|+h_k^2K\mu+(1+h_kK)^2h_k\mu,
\end{split}
\end{eqnarray*}
which justifies the validity of \eqref{e:38} at $j=3$. Suppose next that \eqref{e:38} holds for $t_j$ as $j\ge 3$ and show that it is also satisfied for $t_{j+1}$. Indeed, employing \eqref{e:26} and \eqref{e:36} tells us that
\begin{eqnarray*}
\begin{split}
\left\|x^k(t_{j+1})-\ox(t_{j+1})\right\|=&\left\|x^k(t_{j+1})-y^k_1(t_{j+1})\right\|\\
\le&(1+h_kK)\left\|x^k(t_j)-\ox(t_j)\right\|+h_k\left\|\dfrac{\ox(t_{j+1})-\ox(t_j)}{t_{j+1}-t_j}-\dot\ox(t_j)\right\|\\
\le&(1+h_kK)\bigg(h_k\sum^{j-1}_{i=0}\left\|\dfrac{\ox(t_{i+1})-\ox(t_i)}{t_{i+1}-t_i}-\dot\ox(t_i)\right\|\\
+&h_k^2K\mu\sum_{i=0}^{j-3}(1+h_kK)^i+(1+h_kK)^{j-1}h_k\mu\bigg)+h_k\left\|\dfrac{\ox(t_{j+1})-\ox(t_j)}{t_{j+1}-t_j}-\dot\ox(t_j)\right\|\\
\le&h_k\sum^j_{i=0}\left\|\dfrac{\ox(t_{i+1})-\ox(t_i)}{t_{i+1}-t_i}-\dot\ox(t_i)\right\|+h_k^2K\mu+h_k^2K\mu\sum_{i=1}^{j-2}(1+h_kK)^i+(1+h_kK)^jh_k\mu\\
=&h_k\sum^j_{i=0}\left\|\dfrac{\ox(t_{i+1})-\ox(t_i)}{t_{i+1}-t_i}-\dot\ox(t_i)\right\|+h_k^2K\mu\sum_{i=0}^{j-2}(1+h_kK)^i+(1+h_kK)^jh_k\mu,
\end{split}
\end{eqnarray*}
which shows that estimate \eqref{e:38} holds for $t_{j+1}$, and thus it is justified for all $j=3,\ldots,k$. Now picking any $j\in\{3,\ldots,k\}$ and using the first inequality in \eqref{e:26}, we get
\begin{eqnarray*}
\begin{split}
\left\|x^k(t_j)-\ox(t_j)\right\|\le&h_k\mu+h_k\mu[(1+h_kK)^{j-2}-1]+(1+h_kK)^kh_k\mu\\
\le&2h_k\mu(1+h_kK)^k=2h_k\mu\left(1+\dfrac{K}{k}\right)^k\le 2h_k\mu e^K.
\end{split}
\end{eqnarray*}
Combining it with \eqref{e:37}, we arrive at \eqref{e:39}. This readily implies that $1-\tau-\varepsilon_k\le\|u^k(t_j)\|\le 1+\tau+\varepsilon_k$ for $j\le j_\tau(k)-1$ and $j\ge j^\tau(k)+1$, i.e., the relationships in \eqref{e:28} are satisfied with $\ve_k$ defined in \eqref{e:39}. Furthermore, it follows from \eqref{e:39} and \eqref{e:36} that
\begin{equation}\label{e:40}
\begin{aligned}
\left\|x^k(t)-y^k_1(t)\right\|\le&(1+h_kK)2h_k\mu e^K+h_k\left\|\dfrac{\ox(t_{j+1})-\ox(t_j)}{t_{j+1}-t_j}-\dot\ox(t_i)\right\|\\
\le&2h_k\mu e^K(1+h_kK)+h_k\mu\;\mbox{ for }\;t\in[t_j,t_{j+1}]\;\mbox{ and }\;j=0,\ldots,k-1.
\end{aligned}
\end{equation}\vspace*{-0.1in}
Next we consider relationships for the $u$-component of $z^k(\cdot)$. The first and third conditions in \eqref{e:26} yield
\begin{eqnarray*}
\begin{aligned}
\sum^{k-2}_{j=0}\bigg\|&\dfrac{u^k(t_{j+2})-u^k(t_{j+1})}{t_{j+2}-t_{j+1}}-\dfrac{u^k(t_{j+1})-u^k(t_{j})}{t_{j+1}-t_{j}}\bigg\|\\
\le&\sum^{k-2}_{j=0}\left\|\dfrac{\ou(t_{j+2})-\ou(t_{j+1})}{t_{j+2}-t_{j+1}}-\dfrac{\ou(t_{j+1})-\ou(t_{j})}{t_{j+1}-t_{j}}\right\|
+2\sum^{k-1}_{j=0}\left\|\dfrac{x^k(t_{j+1})-x^k(t_{j})}{t_{j+1}-t_{j}}-\dfrac{\ox(t_{j+1})-\ox(t_{j})}{t_{j+1}-t_{j}}\right\|\\
\le&\mu+2\sum^{k-1}_{j=0}\left\|-v^k_j-\dfrac{\ox(t_{j+1})-u^k(t_{j})}{t_{j+1}-t_{j}}\right\|=\mu+2\sum^{k-1}_{j=0}\left\|-v^k_j-w^k_{1j}\right\|
=\mu+2\sum^{k-1}_{j=0}\dist\left(-\dfrac{\ox(t_{j+1})-\ox(t_j)}{t_{j+1}-t_j};F\big(z^k(t_j)\big)\right)\\
\le&\mu+2\sum^{k-1}_{j=0}\left\|\dfrac{\ox(t_{j+1})-\ox(t_j)}
{t_{j+1}-t_j}-\dot{\ox}(t_j)\right\|+2\sum^{k-1}_{j=0}K\left\|x^k(t_j)-\ox(t_j)\right\|\\
\le&3\mu+2Kk2h_k\mu e^K=3\mu+4K\mu e^K\le\Tilde{\mu}=\max\big\{3\mu+4K\mu e^K,4\mu e^K+\mu\big\},
\end{aligned}
\end{eqnarray*}
which justifies the first estimate in \eqref{e:31}. To verify the second estimate therein, we deduce from \eqref{e:33}, \eqref{e:39}, and the second inequality in \eqref{e:26} that
\begin{eqnarray*}
\begin{aligned}
\left\|\dfrac{u^k(t_1)-u^k(t_0)}{t_1-t_0}\right\|\le&\left\|\dfrac{u^k(t_1)-\ou(t_1)}{t_1-t_0}\right\|+\left\|\dfrac{u^k(t_0)-\ou(t_0)}{t_1-t_0}
\right\|+\left\|\dfrac{\ou(t_1)-\ou(t_0)}{t_1-t_0}\right\|\\
\le&\left\|\dfrac{x^k(t_1)-\ox(t_1)}{t_1-t_0}\right\|+\left\|\dfrac{x^k(t_0)-\ox(t_0)}{t_1-t_0}\right\|+\mu\le 4\mu e^K+\mu\le\Tilde{\mu},
\end{aligned}
\end{eqnarray*}
which readily gives us the claimed result in \eqref{e:31}.\vspace*{-0.05in}

It remains to justify the $W^{1,2}$-convergence of $z^k(t)$ to $\oz(t)$ in \eqref{e:30}. Using \eqref{e:33} for $j=0$ with $x^k(t_0)=x_0$, the construction of $a^k(t)$, and the Newton-Leibniz formula, it suffices to show that the sequence of $(\dot x^k(t),\dot u^k(t))$ converges to $(\dot\ox(t),\dot\ou(t))$ strongly in $L^2[0,T]$. To this end we have
\begin{eqnarray*}
\begin{aligned}
\int^T_0\left\|\dot{x}^k(t)-w^k_1(t)\right\|^2dt=&\sum^{k-1}_{j=0}(t_{j+1}-t_j)\left\|-v^k_j-w^k_{1j}\right\|^2=\sum^{k-1}_{j=0}h_k\dist^2\left(-\dfrac{\ox(t_{j+1})
-\ox(t_j)}{t_{j+1}-t_j};F\big(z^k(t_j)\big)\right)\\
\le&\sum^{k-1}_{j=0}h_k\left(K\left\|x^k(t_j)-\ox(t_j)\right\|+\left\|\dfrac{\ox(t_{j+1})-\ox(t_j)}{t_{j+1}-t_j}-\dot{\ox}(t_j)\right\|\right)^2\\
\le&2\sum^{k-1}_{j=0}h_kK^2\left\|x^k(t_j)-\ox(t_j)\right\|^2+2\sum^{k-1}_{j=0}h_k\left\|\dfrac{\ox(t_{j+1})-\ox(t_j)}{t_{j+1}-t_j}-\dot{\ox}(t_j)
\right\|^2\\\le&2\sum^{k-1}_{j=0}h_kK^2\left(2h_k\mu e^K\right)^2+2\sum^{k-1}_{j=0}h_k\left\|\dfrac{\ox(t_{j+1})-\ox(t_j)}{t_{j+1}-t_j}-\dot{\ox}(t_j)\right\|^2\\
\le&8K^2h_k^2\mu^2e^{2K}+2\sum^{k-1}_{j=0}h_k\left\|\dfrac{\ox(t_{j+1})-\ox(t_j)}{t_{j+1}-t_j}-\dot{\ox}(t_j)\right\|^2\to 0
\end{aligned}
\end{eqnarray*}
as $k\to\infty$ due to \eqref{e:25}, \eqref{yk}, and the definition of $w^k(t)$. It follows furthermore that
\begin{eqnarray*}
\begin{aligned}
\int^T_0\left\|\dot{u}^k(t)-w^k_2(t)\right\|^2dt=&\int^T_0\left\|\dfrac{u^k(t_{j+1})-u^k(t_j)}{h_k}-\dfrac{\ou(t_{j+1})-\ou(t_j)}{h_k}\right\|^2dt\\
=&\int^T_0\left\|\dfrac{u^k(t_{j+1})-\ou(t_{j+1})}{h_k}-\dfrac{u^k(t_{j})-\ou(t_j)}{h_k}\right\|^2dt\\
=&\int^T_0\left\|\dfrac{x^k(t_{j+1})-\ox(t_{j+1})}{h_k}-\dfrac{x^k(t_{j})-\ox(t_j)}{h_k}\right\|^2dt\\
=&\int^T_0\left\|\dfrac{x^k(t_{j+1})-x^k(t_j)}{h_k}-\dfrac{\ox(t_{j+1})-\ox(t_j)}{h_k}\right\|^2dt\\
=&\int^T_0\left\|\dot{x}^k(t)-w^k_1(t)\right\|^2dt\to 0\;\mbox{ as }\;k\to\infty
\end{aligned}
\end{eqnarray*}
due to the above convergence of $\{\dot x^k(\cdot)\}$. This verifies \eqref{e:30} and completes the proof of the theorem. $\h$\vspace*{-0.15in}

\section{Existence of Optimal Sweeping Solutions and Relaxation}
\setcounter{equation}{0}\vspace*{-0.1in}

In this section we start studying optimal solutions to the original sweeping control problem $(P)$. By taking into account the discussion above and further applications, our main attention is paid to the parametric family of {\em problems $(P^\tau)$} as $\tau\ge 0$ with $(P^0)=(P)$, which are different from $(P)$ only in that the control constraint \eqref{e:7} is replaced by those in \eqref{e:8}. First we establish the following existence theorem of optimal solutions for $(P^\tau)$ in the the class of $W^{1,2}[0,T]$ functions.\vspace*{-0.1in}

\begin{theorem} {\bf (existence of sweeping optimal solutions).}\label{Th:2} Given $r>0$ and $T>0$, consider the optimal control problem $(P^{\tau})$ for any fixed $\tau\in[0,\Bar\tau]$ as $\Bar\tau:=\min\{r,T\}$ in the equivalent form of the differential inclusion \eqref{e:21} over all the $W^{1,2}[0,T]$ triples $z(\cdot)=(x(\cdot),u(\cdot),a(\cdot))$. In addition to the assumptions in {\rm(H1)} and {\rm(H2)}, suppose that along some minimizing sequence of $z^k(\cdot)=(x^k(\cdot),u^k(\cdot),a^k(\cdot))$, $k\in\N$, we have that $\{\dot u^k(\cdot)\}$ is bounded in $L^2([0,T];\R^n)$ while $\{a^k(\cdot)\}$ is bounded in $W^{1,2}([0,T];\R^d)$ and that the running cost $\ell$ in \eqref{e:4} is convex with respect to the velocity variables $(\dot x,\dot u,\dot a)$. Then each sweeping control problem $(P^\tau)$ admits an optimal solution.
\end{theorem}\vspace{-0.1in}
{\bf Proof.} Fix any $\tau\in[0,\Bar\tau]$ and deduce from Proposition~\ref{Th:1} that the set of feasible solutions to $(P^\tau)$ is nonempty. It follows from the assumption imposed on $\{(u^k(\cdot),a^k(\cdot))\}$ by basic functional analysis that the sequence $\{(\dot u^k(\cdot),\dot a^k(\cdot))\}$ is weakly compact in $L^2([0,T[;R^{n\times d})$. Thus there are functions $\vartheta^u(\cdot)\in L^2([0,T];\R^n)$ and $\vartheta^a(\cdot)\in L^2([0,T];\R^d)$ such that $\dot{u}^k(\cdot)\to\vartheta^u(\cdot)$ and $\dot{a}^k(\cdot)\to\vartheta^a(\cdot)$ along some subsequence $k\to\infty$ weakly in $L^2([0,T];\R^n)$ and $L^2([0,T];\R^d)$, respectively. By taking into account that $\|u^k(0)\|=r$ by \eqref{e:8} and that the sequence $\{a^k(0)\}$ is bounded, we can assume without loss of generality that $u^k(0)\to u_0$ and $a^k(0)\to a_0$ as $k\to\infty$ for some $u_0\in\R^n$ and $a_0\in\R^d$. Defining now the absolutely continuous functions $\ou:[0,T]\to\R^n$ and $\oa:[0,T]\to\R^d$ by
\begin{equation}\label{e:16}
\ou(t):=u_0+\int^t_0\vartheta^u(s)ds\;\mbox{ and }\;\oa(t):=a_0+\int^t_0\vartheta^a(s)ds,
\end{equation}
we see that $(u^k(\cdot),a^k(\cdot))\to(\ou(\cdot),\oa(\cdot))$ in the norm of $W^{1,2}([0,T];\R^{n\times d})$. This implies that $(\ou(\cdot),\oa(\cdot))\in W^{1,2}([0,T];\R^{n\times d})$ and that $\ou(\cdot)$ satisfies the constraints in \eqref{e:8}. Furthermore, it follows from Proposition~\ref{Th:1} that the trajectories $x^k(\cdot)$ of \eqref{e:21} uniquely generated by $(u^k(\cdot),a^k(\cdot))$ are uniformly bounded in $W^{1,2}([0,T];\R^n)$, and hence a subsequence of them converges to some $\ox(\cdot)\in W^{1,2}([0,T];\R^n)$.\vspace*{-0.05in}

Let us show that the limiting triple $\oz(\cdot)$ satisfies \eqref{e:21} with $F(z)$ defined in \eqref{e:17} and that
\begin{equation}\label{lsc}
J[\ox,\ou,\oa]\le\liminf_{k\to\infty}J[x^k,u^k,a^k]
\end{equation}
for the cost functional \eqref{e:4}. To proceed, we apply the Mazur weak closure theorem to the sequence $\{\dot z^k(\cdot)\}$, which tells us that the sequence of convex combination of $\dot z^k(\cdot)$ converges to $\dot\oz(\cdot)$ weakly in $L^2[0,T]$, and so its subsequence converges to $\dot\oz(t)$ for a.e.\ $t\in[0,T]$. It follows from the above that $\oz(\cdot)$ satisfies the differential inclusion \eqref{e:21} due to the convexity of the sets $F(z)$. Using finally the imposed convexity of the running cost $\ell$ in $\dot z$ and the assumptions in (H2) together with the Lebesgue dominated convergence theorem yields \eqref{lsc} and thus completes the proof of the theorem. $\h$\vspace*{-0.02in}

We can see that the underlying assumption of Theorem~\ref{Th:2} is the {\em convexity} of the integrand $\ell$ with respect to velocities. This assumption, which is not needed for deriving necessary optimality conditions, can be generally relaxed (and even fully dismissed in rather broad nonconvex settings from the viewpoint of actual solving optimization problems for differential inclusions) due the so-called {\em Bogoluybov-Young relaxation procedure}. To describe it in the setting of $(P^\tau)$, denote by $\ell_{F}(t,x,u,a,\dot x,\dot u,\dot a)$ the convexification of the integrand in \eqref{e:4} on the set $F(x,u,a)$ from \eqref{e:17} with respect to the velocity variables $(\dot x,\dot u,\dot a)$ for all $t,x,u,a$, i.e., the largest convex and l.s.c.\ function majorized by $\ell(t,x,u,a,\cdot,\cdot,\cdot)$ on this set; we put $\Hat\ell:=\infty$ at points out of $F(x,u,a)$. Define now the {\em relaxed sweeping problem} $(R^\tau)$ by
\begin{equation}\label{relax}
\mbox{minimize }\;\Hat{J}[z]:=\varphi\big(x(T)\big)+\int^{T}_{0}\Hat{\ell}_{F}\big(t,x(t),u(t),a(t),\dot{x}(t),\dot{u}(t),\dot{a}(t)\big)dt
\end{equation}
over all the triples $z(\cdot)=(x(\cdot),u(\cdot),a(\cdot))\in W^{1,2}[0,T]$ satisfying the constraints in \eqref{e:8}. Of course, there is no difference between problems $(P^\tau)$ and $(R^\tau)$ if the integrand $\ell$ is convex with respect to $(\dot x,\dot u,\dot a)$. Furthermore, Theorem~\ref{Th:2} ensures the existence of optimal solutions to $(R^\tau)$. The strong relationship between the original and relax/convexified problems, known as {\em relaxation stability}, is that in many situations the optimal values of the cost functionals therein agree. This phenomenon has been well recognized for differential inclusions with Lipschitzian right-hand sides in state variables (see \cite{t}), which is never the case for the sweeping process. A more subtle result of this type is obtained in \cite[Theorem~4.2]{dfm} for differential inclusions satisfying the modified one-sided Lipschitz property, which however is also restrictive in applications to sweeping control. The relaxation stability result that directly concerns sweeping control problems is given in \cite[Theorem~2]{et} while it deals only with the case of controlled perturbations. In general, relaxation stability in sweeping optimal control is an open question.\vspace*{-0.05in}

Our current study here and its continuation in \cite{cm1} concern local optimal solutions to $(P^\tau)$ involving a local version of relaxation stability. Following \cite{m95}, we say that $\bar{z}(\cdot)$ is a {\em relaxed intermediate local minimizer} (r.i.l.m) for $(P^\tau)$ if it is feasible to this problem with $J[\bar{z}]=\Hat{J}[\bar{z}]$ and if there are numbers $\alpha\ge 0$ and $\epsilon>0$ such that $J[\bar{z}]\le J[z]$ for any feasible solution $z(\cdot)$ to $(P^\tau)$ satisfying
\begin{equation}\label{ilm}
\left\|z(t)-\oz(t)\right\|<\epsilon\;\mbox{ for all }\;t\in[0,T]\;\mbox{ and }\;\alpha\int^{T}_{0}\left\|\dot{z}(t)-\dot{\oz}(t)\right\|^{2}dt<\epsilon.
\end{equation}
This notion distinguishes local minimizers that lie between classical weak and strong minima in continuous-time variational problems and can be strictly different from both of them even in fully convex settings; see \cite{m-book2} for discussions, examples, and references. It is clear that from the viewpoint of deriving necessary optimality conditions we can confine ourselves to the case of $\alpha=1$.\vspace*{-0.15in}

\section{Discrete Approximations of Local Optimal Solutions}
\setcounter{equation}{0}\vspace*{-0.1in}

In this section we construct a sequence of well-posed discrete approximations of each problem $(P^\tau)$ as $0\le\tau\le\Bar\tau$ with $\Bar\tau=\min\{r,T\}$ and then employ this method to the study of relaxed intermediate local minimizers for this problem. Given any r.i.l.m.\ $\bar{z}=(\bar{x}(\cdot),\bar{u}(\cdot),\bar{a}(\cdot))$ for $(P^\tau)$ and the discrete mesh $\Delta_k$ from \eqref{mesh}, for every $k\in\N$ define the {\em discrete sweeping control problem} $(P^\tau_k)$ as follows: minimize
\begin{eqnarray}\label{disc-cost}
\begin{aligned}
\disp J_{k}[z^{k}]:=&\varphi(x^{k}_{k})+h_{k}\sum^{k-1}_{j=0}\ell\left(t^{k}_{j},x^{k}_{j},u^{k}_{j},a^{k}_{j},\dfrac{x^{k}_{j+1}-x^{k}_{j}}{h_{k}},\dfrac{u^{k}_{j+1}
\disp-u^{k}_{j}}{h_{k}},\dfrac{a^{k}_{j+1}-a^{k}_{j}}{h_{k}}\right)+\Big\|\frac{x^k_1-x^k_0}{h_k}-\dot\ox(0)\Big\|^2\\
\disp+&\sum^{k-1}_{j=0}\int^{t^{k}_{j+1}}_{t^{k}_{j}}{\left(\left\|\dfrac{x^{k}_{j+1}-x^{k}_{j}}{h_{00k}}-\dot\ox(t)\right\|^{2}+\left\|\dfrac{u^{k}_{j+1}
\disp-u^{k}_{j}}{h_{k}}-\dot{\bar{u}}(t)\right\|^{2}+\left\|\dfrac{a^{k}_{j+1}-a^{k}_{j}}{h_{k}}-\dot{\bar{a}}(t)\right\|^{2}\right)}dt\\
\disp+&\dist^2\left(\left\|\dfrac{u^k_1-u^k_0}{h_k}\right\|;\big(-\infty,\Tilde{\mu}\big]\right)+\dist^2\left(\sum^{k-2}_{j=0}\left\|\dfrac{u^k_{j+2}-2u^k_{j+1}
\disp+u^k_j}{h_k}\right\|;\big(-\infty,\Tilde{\mu}\big]\right)
\end{aligned}
\end{eqnarray}
over elements $z^{k}:=(x^{k}_{0},x^{k}_{1},\ldots,x^{k}_{k},u^{k}_{0},u^{k}_{1},\ldots,u^{k}_{k-1},a^{k}_{0},a^{k}_{1},\ldots,a^{k}_{k-1})$ satisfying the constraints
\begin{equation}\label{e:41}
x^{k}_{j+1}\in x^{k}_{j}-h_{k}F(x^{k}_{j},u^{k}_{j},a^{k}_{j})\;\mbox{ for }\;j=0,\ldots,k-1\;\mbox{ with }\;(x^{k}_{0},u^k_0,a^k_0)=\big(x_{0},\ou(0),\oa(0)\big),
\end{equation}
\begin{equation}\label{e:42}
\left\langle x^*_{i},x^{k}_k-u^{k}_k\right\rangle\le 0\;\mbox{ for }\;i=1,\ldots,m,
\end{equation}
\begin{equation}\label{e:43}
\|u^k_j\|=r\;\mbox{for}\;j=j_\tau(k),\ldots,j^\tau(k);\;r-\tau-\varepsilon_k\le\|u^k_j\|\le r+\tau+\varepsilon_k\;\mbox{for}\;j\le j_\tau(k)-1\;\mbox{and}\;j\ge j^\tau(k)+1,
\end{equation}
\begin{equation}\label{e:44}
\left\|(x^{k}_{j},u^{k}_{j},a^{k}_{j})-\big(\bar{x}(t^{k}_{j}),\bar{u}(t^{k}_{j}),\bar{a}(t^{k}_{j})\big)\right\|\le\epsilon/2\;\mbox{ for }\;j=0,\ldots,k-1,
\end{equation}
\begin{equation}\label{e:45}
\sum^{k-1}_{j=0}\int^{t^{k}_{j+1}}_{t^{k}_{j}}{\left(\left\|\dfrac{x^{k}_{j+1}-x^{k}_{j}}{h_{k}}-\dot{\bar{x}}(t)\right\|^{2}+\left\|\dfrac{u^{k}_{j+1}-u^{k}_{j}}
{h_{k}}-\dot{\bar{u}}(t)\right\|^{2}+\left\|\dfrac{a^{k}_{j+1}-a^{k}_{j}}{h_{k}}-\dot{\bar{a}}(t)\right\|^{2}\right)}dt\le\dfrac{\epsilon}{2},
\end{equation}
\begin{equation}\label{e:46}
\left\|\dfrac{u^k_1-u^k_0}{t^k_1-t^k_0}\right\|\le\Tilde{\mu}+1,\;\mbox{ and }\;\sum^{k-2}_{j=0}\left\|\dfrac{u^k_{j+2}-2u^k_{j+1}+u^k_j}{h_k}\right\|\le \Tilde{\mu}+1,
\end{equation}
where $\epsilon>0$ is taken from definition \eqref{ilm} with $\alpha=1$ while $\ve_k$ and $\Tilde\mu$ are taken from Theorem~\ref{Th:3}.\vspace*{-0.05in}

Let us first show that each problem $(P^\tau_k)$ admits an optimal solution for all large $k\in\N$; this issue is unavoidable in employing the method of discrete approximations to study local minimizers for $(P^\tau)$.\vspace*{-0.1in}

\begin{proposition}{\bf(existence of optimal solutions to discrete approximations).}\label{Th:4} Suppose that {\rm(H1)} holds and that {\rm(H2)} is also satisfied around the given local minimizer $\oz(\cdot)$ for $(P^\tau)$. Then each problem $(P^\tau_k)$ admits an optimal solution provided that $k\in\N$ is sufficiently large.
\end{proposition}\vspace*{-0.1in}
{\bf Proof.} Theorem~\ref{Th:3} tells us that the set of feasible solutions to $(P^\tau_k)$ is nonempty for all large $k\in\N$. Moreover, the constraints in \eqref{e:43}--\eqref{e:45} ensure that this set is bounded. To justify the claimed existence of optimal solutions to $(P^\tau_k)$ by the Weierstrass existence theorem, it remains to verify that this set is closed. To proceed, take a sequence $z^\nu(\cdot)=z^{\nu}:=(x^{\nu}_{0},\ldots,x^{\nu}_{k},u^{\nu}_{0},\ldots,u^{\nu}_{k-1},a^{\nu}_{0},\ldots,a^{\nu}_{k-1})$ of feasible solutions for $(P^\tau_k)$ converging to some $z(\cdot)=z:=(x_{0},\ldots,x_{k},u_{0},\ldots,u_{k-1},a_{0},\ldots,a_{k-1})$ as $\nu\to\infty$ and show that $z$ is feasible to $(P^\tau_k)$ as well. Observe that $\left\langle x^*_i,x_{j}-u_{j}\right\rangle=\lim_{\nu\to\infty}\left\langle x^*_i,x^{\nu}_{j}-u^{\nu}_{j}\right\rangle\le 0$ for all $i=1,\ldots,m,$ and $j=0,\ldots,k-1$, and so $x_{j}-u_{j}\in C$ for all $j=0,\ldots,k-1$. Picking now
$i\in\{1,\ldots,m\}\backslash I(x_{j}-u_{j})$, we have $\left\langle x^*_i,x_{j}-u_{j}\right\rangle<0$, which yields $\left\la x^*_i,x^{\nu}_{j}-u^{\nu}_{j}\right\ra<0$ for $\nu$ sufficiently large. Then it follows that $i\in\{1,\ldots,m\}\backslash I(x^{\nu}_{j}-u^{\nu}_{j})$ and hence $I(x^{\nu}_{j}-u^{\nu}_{j})\subset I(x_{j}-u_{j})$ for $\nu\in\N$ sufficiently large. By taking \eqref{e:17} and \eqref{e:20} into account, we get the equalities
$$
x^{\nu}_{j+1}-x^{\nu}_{j}=-h_{k}\left(\sum_{i\in I(x^{\nu}_{j}-u^{\nu}_{j})}\lambda^{\nu}_{ji}x^{*}_{i}+f(x^{\nu}_{j},a^{\nu}_{j})\right)
=-h_{k}\left(\sum_{i\in I(x_{j}-u_{j})}\lambda^{\nu}_{ji}x^*_{i}+f(x^{\nu}_{j},a^{\nu}_{j})\right),
$$
where $\lambda^{\nu}_{ji}:=0$ if $i\in I(x_{j}-u_{j})\backslash I(x^{\nu}_{j}-u^{\nu}_{j})$. This shows therefore that
\begin{equation*}
\dfrac{x^{\nu}_{j+1}-x^{\nu}_{j}}{-h_{k}}-f(x^{\nu}_{j},a^{\nu}_{j})=\sum_{i\in I(x_{j}-u_{j})}\lambda^{\nu}_{j}x^{*}_{i}\in N(x_{j}-u_{j};C).
\end{equation*}
Passing there to the limit as $\nu\to\infty$ and using the closedness of $N(x_{j}-u_{j};C)$ give us
$$
\dfrac{x_{j+1}-x_{j}}{-h_{k}}-f(x_{j},a_{j})\in N(x_{j}-u_{j};C),
$$
which ensures that $x_{j+1}\in x_{j}-h_{k}F(x_{j},u_{j},a_{j})$ and thus completes the proof of the proposition. $\h$

The next theorem is a key result of the method of discrete approximations in sweeping optimal control. It shows that optimal solutions to $(P^\tau)$ and $(P^\tau_k)$ are so closely related that solving the continuous-time control problem $(P^\tau)$ for small $\tau\ge 0$ can be practically replaced by solving its finite-dimensional discrete counterparts $(P^\tau_k)$ when $k$ is sufficiently large. Moreover, it justifies the possibility to derive necessary optimality conditions for local minimizers of $(P^\tau)$ by passing to the limit from those in $(P^\tau_k)$ as $k\to\infty$.\vspace*{-0.05in}

\begin{theorem}{\bf(strong discrete approximation of intermediate local minimizers).}\label{Th:5} Let $\bar{z}(\cdot)=(\bar{x}(\cdot),\bar{u}(\cdot),\bar{a}(\cdot))$ be a r.i.l.m.\ for problem $(P^\tau)$, where $\tau\in[0,\Bar\tau]$ with $\Bar\tau=\min\{r,T\}$. In addition to the assumptions in Theorem~{\rm\ref{Th:3}} and Proposition~{\rm\ref{Th:4}} imposed on $\oz(\cdot)$, suppose that both terminal and running costs in \eqref{e:4} are continuous at $\bar{x}(T)$ and at $(t,\bar{z}(t),\dot{\bar{z}}(t))$ for a.e.\ $t\in[0,T]$, respectively, and that $\ell(\cdot,z,\dot z)$ is uniformly majorized by a summable function near the given local minimizer. Then any sequence of piecewise linearly extended to $[0,T]$ optimal solutions $\oz^k(\cdot)=(\ox^k(\cdot),\ou^k(\cdot),\oa^k(\cdot))$ of $(P^\tau_k)$ converges to $\bar{z}(\cdot)$ in the norm topology of $W^{1,2}([0,T];\R^{2n+d})$ with
\begin{equation}\label{x01}
\frac{\ox^k_1-\ox^k_0}{h_k}\to\dot\ox(0)\;\mbox{ as }\;k\to\infty
\end{equation}
and the validity of the estimates
\begin{equation}\label{e:50}
\left\|\dfrac{\ou^k_1-\ou^k_0}{h_k}\right\|\le\Tilde\mu,\quad\limsup_{k\to\infty}\sum^{k-2}_{j=0}\left\|\dfrac{\ou^k_{j+2}-2\ou^k_{j+1}+\ou^k_j}{h_k}\right\|\le\Tilde{\mu},
\end{equation}
where the number $\Tilde\mu$ is calculated in \eqref{e:31}.
\end{theorem}\vspace*{-0.1in}
{\bf Proof.} Fix a sequence of optimal solutions $\oz^k(\cdot)$ to $(P^\tau_k)$, which exists by Proposition~\ref{Th:4}. It is easy to see that all the statements of the theorem are implied by the equality
\begin{equation}\label{e:52}
\begin{aligned}
\lim_{k\to\infty}&\int^{T}_{0}\left(\left\|\dot{\ox}(t)-\dot{{\ox}}^{k}(t)\right\|^{2}+\left\|\dot{\ou}(t)-\dot{{\ou}}^{k}(t)\right\|^{2}+
\left\|\dot{\oa}(t)-\dot{{\oa}}^{k}(t)\right\|^{2}\right)dt+\Big\|\frac{\ox^k_1-\ox^k_0}{h_k}-\dot\ox(0)\Big\|^2\\
+&\dist^2\left(\left\|\dfrac{\ou^k_1-\ou^k_0}{h_k}\right\|;\big(-\infty,\Tilde{\mu}\big]\right)+\dist^2\left(\sum^{k-2}_{j=0}\left\|\dfrac{\ou^k_{j+2}-2\ou^k_{j+1}
+\ou^k_j}{h_k}\right\|;\big(-\infty,\Tilde{\mu}\big]\right)=0.
\end{aligned}
\end{equation}
To justify \eqref{e:52}, suppose that it does not hold, i.e., there is a subsequence of natural numbers (without relabeling) along which the limit in \eqref{e:52} equals to some $c>0$. By the weak compactness of the unit ball in $L^{2}([0,T];\R^{2n+d})$ we can find a triple $(v(\cdot),w(\cdot),q(\cdot))\in L^{2}([0,T];\R^{2n+d})$ and yet another subsequence of $\{z^k(\cdot)\}$--again without relabeling--such that
$$
\big(\dot{{\ox}}^{k}(\cdot),\dot{{\ou}}^{k}(\cdot),\dot{{\oa}}^{k}(\cdot)\big)\to\big(v(\cdot),w(\cdot),q(\cdot)\big)\;\mbox{ weakly in }\;L^{2}([0,T];\R^{2n+d}).
$$
Define now the absolutely continuous function $\tz(\cdot):=(\tx(\cdot),\tu(\cdot),\ta(\cdot))\colon[0,T]\to\R^{2n+d}$ by
$$
\tz(t):=\big(x_{0},\bar{u}(0),\bar{a}(0)\big)+\int^{t}_{0}(v(s),w(s),q(s)\big)ds,\quad t\in[0,T],
$$
which gives us $\dot{\tz}(t)=(v(t),w(t),q(t))$ a.e.\ on $[0,T]$ and implies that $\dot\oz^k(\cdot)\to\dot{\tz}(\cdot)=(\dot{\tx}(\cdot),\dot{\tu}(\cdot),\dot{\ta}(\cdot))$ weakly in $L^{2}([0,T];\R^{2n+d})$ and therefore ${\oz}^k(\cdot)\to\tz(\cdot)\in W^{1,2}([0,T];\R^{2n+d})$ in the $W^{1,2}$-norm topology. Furthermore, it follows from the Mazur weak closure theorem that there is a sequence of convex combinations of $\dot{{\oz}}^{k}(\cdot)$ that converges to $\dot{\tz}(\cdot)$ strongly in $L^{2}([0,T];\R^{2n+d})$ and thus almost everywhere on $[0,T]$ along a subsequence. It is clear that the limiting $u$-component $\tu(\cdot)$ obeys the constraints in \eqref{e:8}. Let us verify that $\tz(\cdot)$ satisfies the differential inclusion \eqref{e:5}, where the moving set $C(t)$ is generated by $\tu(\cdot)$ in \eqref{e:6}.\vspace*{-0.05in}

To proceed, observe first that $\tx(t)-\tu(t)=\lim_{k\to\infty}({\ox}^{k}(t)-{\ou}^{k}(t))\in C$ by the closedness of the polyhedron $C$. It follows from the above that there are a function $\nu:\N\to\N$ and a sequence of real numbers $\{\alpha(k)_j|\;j=k,\ldots,\nu(k)\}$ such that
$$
\alpha(k)_j\ge 0,\;\sum_{j=k}^{\nu(k)}\alpha(k)_j=1,\;\mbox{ and }\;\sum_{j=k}^{\nu(k)}\alpha(k)_j\dot{\tz}^j(t)\to\dot{\tz}(t)\;\mbox{ a.e. }\;t\in [0,T]
$$
as $k\to\infty$. Then by the closedness and convexity of the normal cone we have the relationships
\begin{equation*}
\begin{aligned}
-\dot{\tx}(t)-f\big(\tx(t),\ta(t)\big)=&\lim_{k\to\infty}\bigg(-\sum_{j=k}^{\nu(k)}\alpha(k)_j\dot{\ox}^j(t)-\sum_{j=k}^{\nu(k)}\alpha(k)_j f\big(\ox^j(t),\oa^j(t)\big)\bigg)\\
=&\lim_{k\to\infty}\sum_{i\in I(\tx(t)-\tu(t))}\bigg(\sum_{j=k}^{\nu(k)}\alpha(k)_j\lambda^j_i\bigg)x^*_i\in N\big(\tx(t)-\tu(t);C\big)\;\mbox{ a.e. }\;t\in[0,T],
\end{aligned}
\end{equation*}
where $I(\cdot)$ is taken from \eqref{e:19}, and where $\lambda^{j}_{i}=0$ if $i\in I(\tx(t)-\tu(t))\backslash I(x^j(t)-u^j(t))$ for $j=k,\ldots,\nu(k)$ and all large $k\in\N$. It shows by \eqref{e:20} that $\tz(\cdot)$ satisfies \eqref{e:21} and hence the constraints in \eqref{mixed}.\vspace*{-0.05in}

Consider further the integral functional
$$
{\cal I}[y]:=\int^{T}_{0}\left\|y(t)-\dot{\bar{z}}(t)\right\|^{2}dt
$$
is l.s.c.\ in the weak topology of $L^{2}([0,T];\R^{2n+d})$ due to the convexity of the integrand in $y$. Hence
\begin{equation}\label{e:55}
{\cal I}(\dot{\Tilde z})=\int^{T}_{0}\left\|\dot{\Tilde{z}}(t)-\dot{\bar{z}}(t)\right\|^{2}dt\le\liminf_{k\to\infty}\sum^{k-1}_{j=0}\int^{t^{k}_{j+1}}_{t^{k}_{j}}
\left\|\dfrac{\bar{z}^{k}_{j+1}-\bar{z}^{k}_{j}}{h_{k}}-\dot{\bar{z}}(t)\right\|^{2}dt
\end{equation}
by the construction of $\Tilde z(\cdot)$. Passing to the limit in \eqref{e:44} and \eqref{e:45} as $k\to\infty$ and using \eqref{e:55}, we get
$$
\left\|\tz(t)-\bar{z}(t)\right\|\le\epsilon/2\;\mbox{ on }\;[0,T]\;\mbox{ and }\;\int^{T}_{0}\left\|\dot{\tz}(t)-\dot{\bar{z}}(t)\right\|^{2}dt\le\epsilon/2.
$$
This means that $\tz(\cdot)$ belongs to the given neighborhood of $\bar{z}(\cdot)$ in $W^{1,2}([0,T];\R^{2n+d})$. Furthermore, the definition of $\ell_{F}$ in \eqref{relax} and its convexity in the velocity variables yield
$$
\int^{T}_{0}\Hat{\ell}_{F}\big(t,\tx(t),\tu(t),\ta(t),\dot{\tx}(t),\dot{\tu}(t),\dot{\ta}(t)\big)dt\le\liminf_{k\to\infty}h_{k}\sum^{k-1}_{j=0}\ell \left(t^{k}_{j},\bar{x}^{k}_{j},\bar{u}^{k}_{j},\bar{a}^{k}_{j},\dfrac{\bar{x}^{k}_{j+1}-\bar{x}^{k}_{j}}{h_{k}},\dfrac{\bar{u}^{k}_{j+1}-\bar{u}^{k}_{j}}{h_{k}}.
\dfrac{\bar{a}^{k}_{j+1}-\bar{a}^{k}_{j}}{h_{k}}\right).
$$
Thus the passage to the limit in the cost functional of $(P^\tau_k)$ and the assumption on $c>0$ in the negation of \eqref{e:52} together with (H2) bring us to the relationships
\begin{equation}\label{e:56}
\Hat{J}[\tz]+c=\varphi\big(\tx(T)\big)+\int^{T}_{0}\Hat{\ell}_{F}\big(t,\tx(t),\tu(t),\ta(t),\dot{\tx}(t),\dot{\tu}(t),\dot{\ta}(t)\big)dt+c\le
\liminf_{k\to\infty}J_{k}[\bar{z}^{k}].
\end{equation}
Applying now Theorem~\ref{Th:3} to the r.i.l.m.\ $\oz(\cdot)$ gives us a sequence $\{z^{k}(\cdot)\}$ of feasible solutions to $(P^\tau_k)$ whose extensions to the whole interval $[0,T]$ strongly approximate $\bar{z}(\cdot)$ in the $W^{1,2}$ topology with the additional convergence in \eqref{x0}. Since $\bar{z}^{k}(\cdot)$ is an optimal solution to $(P^\tau_k)$, we have
\begin{equation}\label{e:57}
J_{k}[\bar{z}^{k}]\le J_{k}[z^{k}]\;\mbox{ for each }\;k\in\N.
\end{equation}
It follows from the structure of the cost functionals in $(P^\tau_k)$, the strong $W^{1,2}$-convergence of $z^k(\cdot)\to\oz(\cdot)$ together with \eqref{x0} in Theorem~\ref{Th:3}, and the continuity assumptions on $\varphi$ and $\ell$ imposed in this theorem that $J_{k}[z^{k}]\to J[\bar{z}]$ as $k\to\infty$. Then passing to the limit in \eqref{e:57} gives us
\begin{equation}\label{e:58}
\limsup_{k\to\infty}J_{k}[\bar{z}^{k}]\le J[\bar{z}].
\end{equation}
Combining finally \eqref{e:56} and \eqref{e:58} with $c>0$ and the definition of r.i.l.m., we get
$$
\Hat{J}[\tz]+c\le J[\bar{z}]=\Hat{J}[\bar{z}],\;\mbox{ and so }\;\Hat{J}[\tz]<\Hat{J}[\bar{z}],
$$
which clearly contradicts the fact that $\bar{z}(\cdot)$ is a r.i.l.m.\ for problem $(P^\tau)$. This justifies the validity of \eqref{e:52} and thus completes the proof of the theorem. $\h$\vspace*{-0.15in}

\section{Generalized Differentiation and Calculations}
\setcounter{equation}{0}\vspace*{-0.1in}

After establishing close connections between optimal solutions to the original and discretized sweeping control problems, our further goal is to derive effective necessary optimality conditions to each problem $(P^\tau_k)$ defined in \eqref{disc-cost}--\eqref{e:46}. Looking at this problem, we can see that it is {\em intrinsically nonsmooth}, even if both terminal and running costs are assumed to be  differentiable, which is not the case here. The main reason for this is the unavoidable presence of the geometric constraints \eqref{e:41} with $F$ given by \eqref{e:17} whose increasing number comes from the discretization of the sweeping differential inclusion \eqref{e:5}. We can deal with such problems by using the robust generalized differential constructions, which are basic in variational analysis and its applications; see, e.g., the books \cite{bz,m-book1,rw} and the references therein. Here we first recall their definitions with a brief overview of the needed properties and then deduce from \cite{hmn} major coderivative calculations for the mapping $F$ in \eqref{e:17} via the initial data of the sweeping process. This together with available calculus rules of generalized differentiation plays a crucial role in deriving verifiable necessary optimality conditions for the sweeping control problems under consideration.\vspace*{-0.05in}

Given a set-valued mapping/multifunction $G\colon\R^n\tto\R^m$, denote by
\begin{equation}\label{e:9}
\Limsup_{x\to\ox}G(x):=\big\{y\in\R^m\big|\;\exists\;\mbox{ sequences }\;x_k\to\ox,\;y_k\to y\;\mbox{ such that }\;y_k\in G(x_k),\;k\in\N\big\}
\end{equation}
the (Kuratowski-Painlev\'e) {\em outer limit} of $G$ at $\ox$ with $G(\ox)\ne\emp$. Considering now a set $\Omega\subset\R^n$ locally closed around $\ox\in\Omega$, the (Mordukhovich basic/limiting) {\em normal cone} to $\Omega$ at $\ox$ is defined by
\begin{equation}\label{e:59}
N(\ox;\Omega):= N_{\Omega}(\ox):=\Limsup_{x\to\ox}\left\{\cone[x-\Pi(x;\Omega)]\right\}
\end{equation}
via the outer limit \eqref{e:9}, where $\Pi(x;\Omega)$ stands for the Euclidean projection of $x$ onto $\Omega$. When $\Omega$ is convex, \eqref{e:59} reduces to the normal cone of convex analysis, but it is often nonconvex in nonconvex settings. The crucial feature of \eqref{e:59} and the associated subdifferential and coderivative constructions for functions and multifunctions (see below) is {\em full calculus} based on variational and extremal principles.\vspace*{-0.05in}

For a set-valued mapping $F:\R^n\tto\R^m$ with its graph
\begin{eqnarray*}
\gph F:=\big\{(x,y)\in\R^n\times\R^m\big|\;y\in F(x)\big\}
\end{eqnarray*}
locally closed around $(\ox,\oy)$, the {\em coderivative} of $F$ at $(\ox,\oy)$ generated by \eqref{e:59} is defined by
\begin{equation}\label{cod}
D^*F(\ox,\oy)(u):=\big\{v\in\R^n\big|\;(v,-u)\in N\big((\ox,\oy);\gph F\big)\big\},\;u\in\R^m.
\end{equation}
When $F:\R^n\to\R^m$ is single-valued and continuously differentiable $({\cal C}^1)$ around $\ox$, we have
\begin{eqnarray*}
D^*F(\ox)(u)=\big\{\nabla F(\ox)^*u\big\}\;\mbox{ for all }\;u\in\R^m
\end{eqnarray*}
via the adjoint/transposed Jacobian matrix $\nabla F(\ox)^*$, where $\oy=F(\ox)$ is omitted.\vspace*{-0.05in}

Given an extended-real-valued l.s.c.\ function $\varphi:\R^n\to\Bar{\R}$ with its domain and epigraph
\begin{eqnarray*}
\dom\varphi:=\big\{x\in\R^n\big|\;\varphi(x)<\infty\big\}\;\mbox{ and }\;\epi\varphi:=\big\{(x,\al)\in\R^{n+1}\big|\;\al\ge\varphi(x)\big\},
\end{eqnarray*}
the (first-order) {\em subdifferential} of $\varphi$ at $\ox\in\dom\varphi$ is generated geometrically by \eqref{e:59} as
\begin{equation*}
\partial\varphi(\ox):=\big\{v\in\R^m\;|\;(v,-1)\in N\big((\ox,\varphi(\ox));\epi\varphi\big)\big\}
\end{equation*}
while admitting various equivalent analytic representations that can be found, e.g., in the books \cite{m-book1,rw}.\vspace*{-0.05in}

Our main emphases here is on evaluating the coderivative of the set-valued mapping $F$ from \eqref{e:17} entirely via the given data of the perturbed sweeping process. Note that the partial normal cone structure of the mapping $F$ reveals the {\em second-order} subdifferential nature of the aforementioned construction in the sense of \cite{m92}. For simplicity in further applications, suppose below that the perturbation function $f$ is smooth while observing that the available calculus rules allow us to consider Lipschitzian perturbations.\vspace*{-0.05in}

Having in mind representation \eqref{e:21} of the mapping $F$ in terms of the generating vectors $x^*_i$ of the convex polyhedron \eqref{e:6} with the active constraint indices $I(\ox)$ in \eqref{e:20}, consider the following subsets:
\begin{equation}\label{e:64}
I_0(y):=\big\{i\in I(\bar{x})\big|\;\left\langle x^*_i,y\right\rangle=0\big\}\;\mbox{ and }\;I_>(y):=\big\{i\in I(\bar{x})\big|\;\left\langle x^*_i,y\right\rangle >0\big\},\quad y\in\R^n.
\end{equation}\vspace*{-0.2in}

The next theorem provides an effective upper estimate of the coderivative of $F$ with ensuring the equality therein under an additional assumption on $x^*_i$.\vspace*{-0.1in}

\begin{theorem}{\bf(calculating the coderivative of the sweeping control mapping).}\label{Th:9} Given $F$ in \eqref{e:17} with $C$ from \eqref{e:6}, assume that $f$ is smooth and denote $G(x):=N(x;C)$. Then for any $(x,u,a)\in\R^n\times\R^n\times\R^d$ and $w\in G(x-u)$ we have the coderivative upper estimate
\begin{equation}\label{e:70}
\begin{aligned}
D^*F(x,u,a,w)(y)\subset\bigg\{&z^*\in\R^n\times\R^n\times\R^d\bigg|\;z^*=\bigg(\nabla_x{f}(x,a)^*y+\sum_{i\in I_0(y)\cup I_>(y)}\gamma_{i}x^*_i,\\
-&\sum_{i\in I_0(y)\cup I_>(y)}\gamma_i x^*_{i},\nabla_a f(x,a)^*y\bigg)\bigg\},\quad y\in\dom D^*G\big(x-u,w-f(x,a)\big),
\end{aligned}
\end{equation}
where $I_0(y)$ and $I_>(y)$ are defined in \eqref{e:64} with $\ox=x-u$, and where $\gg_i\in\R$ for $i\in I_0(y)$ while $\gg_i\ge 0$ for $i\in I_>(y)$. Furthermore, \eqref{e:70} holds as an equality provided that the generating vectors $\{x^*_i|\;i\in I(x-u)\}$ of the polyhedron $C$ are linearly independent.
\end{theorem}\vspace*{-0.1in}
{\bf Proof.} Picking any $y\in\dom D^*G(x-u,w-f(x,a))$ and $z^*\in D^*F(x,u,a,X)(y)$ and then denoting $\Tilde{G}(x,u,a):=G(x-u)$ and $\Tilde{f}(x,u,a):=f(x,a)$, we deduce from \cite[Theorem~1.62]{m-book2} that
$$
z^*\in\nabla\Tilde{f}(x,u,a)^*y+D^*\Tilde{G}\big(x-u,w-f(x,a)\big)(y).
$$
Observe then the obvious representation
$$
\Tilde{G}(x,u,a)=G\circ g(x,u,a)\;\mbox{ with }\;g(x,u,a):=x-u,
$$
where the latter mapping has the surjective derivative. It follows from \cite[Theorem~1.66]{m-book2} that
\begin{equation}\label{e:71}
z^*\in\nabla\Tilde{f}(x,u,a)^*y+\nabla g(x,u,a)^*D^*G\big(x-u,w-f(x,a)\big)(y).
\end{equation}
Employing now in \eqref{e:71} the coderivative estimate for the normal cone mapping $G$ obtained in \cite[Theorem~4.5]{hmn} with the exact coderivative calculation given in \cite[Theorem~4.6]{hmn} under the linear independence of the generating vectors $x^*_i$ and also taking into account the structure of the mapping $\Tilde f$ in \eqref{e:71}, we arrive at \eqref{e:70} and the equality therein under the aforementioned assumption. $\h$\vspace*{-0.15in}

\section{Necessary Optimality Conditions}
\setcounter{equation}{0}\vspace*{-0.1in}

In this section we derive necessary conditions for optimal solutions to each discrete approximation problems $(P^\tau_k)$ with $k\in\N$ and $0\le\tau\le\Bar\tau=\min\{r,T\}$. As shown in Theorem~\ref{Th:5}, for large $k\in\N$ and any fixed $\tau\in[0,\Bar\tau]$ the constructed optimal solutions $\oz^k(\cdot)$ to $(P^\tau_k)$ are practically undistinguished (in the $W^{1,2}$ norm) from the optimal solution $\oz(\cdot)$ to the continuous-time sweeping control problem $(P^\tau)$, and so the necessary optimality conditions for $\oz^k(\cdot)$ obtained below can be well treated as ``almost optimality"  necessary conditions for the solution $\oz(\cdot)$ to $(P^\tau)$ playing virtually the same role in applications.\vspace*{-0.05in}

To proceed, we first establish necessary optimality conditions for $(P^\tau_k)$ in the discrete {\em Euler-Lagrange form} via the generalized differential constructions of Section~6; cf.\ \cite{m-book2}. Then employing the complete coderivative calculations of Theorem~\ref{Th:9} for the underlying mapping $F$ from \eqref{e:21} allows us to derive necessary optimality conditions for the sweeping control problem $(P^\tau_k)$ entirely in terms of its initial data. Throughout this section we assume that the cost functions $\varphi$ and $\ell$ in \eqref{disc-cost} are {\em locally Lipschitzian} around the points in question and for brevity drop indicating the time-dependence of the running cost $\ell$.\vspace*{-0.1in}

\begin{theorem}{\bf (Euler-Lagrange conditions for discrete approximations).}\label{Th:10} Fixing any $\tau\in[0,\Bar\tau]$ and $k\in\N$, consider an optimal solution $\oz^{k}=(x_{0},\ox^{k}_{1},\ldots,\ox^{k}_{k},\ou^{k}_{0},\ldots,\ou^{k}_{k},\oa^{k}_{0},\ldots,\oa^{k}_{k})$ to problem $(P^\tau_k)$. Then there exist dual elements $\lambda^{k}\ge0,\;\alpha^k=(\al^k_1,\ldots,\al^k_m)\in\R^m_{+},\;\xi^{k}=(\xi^{k}_{0},\ldots,\xi^{k}_{k})\in\R^{k+1}$, and $p^k_j=(p^{xk}_j,p^{uk}_j,p^{ak}_j)\in\R^n\times\R^n\times\R^d$ as $j=0,\ldots,k$ satisfying the conditions
\begin{equation}\label{e:72}
\lambda^{k}+\left\|\alpha^{k}\right\|+\left\|\xi^{k}\right\|+\sum^{k-1}_{j=0}\left\|p^{xk}_{j}\right\|+\left\|p^{uk}_0\right\|+\left\|p^{ak}_0\right\|\not=0,
\end{equation}
\begin{equation}\label{e:73}
\alpha^{k}_i\left\la x^*_i,\ox^k_k-\ou^k_k\right\ra=0,\quad i=1,\ldots,m,
\end{equation}
\begin{equation}\label{e:74}
\left\{\begin{array}{ll}
\xi^k_j(\nu-\|\ou^k_j\|)\le 0\;\mbox{ for all }\;\nu\in\big[r-\tau-\ve_k,r+\tau+\ve_k\big]\\
\mbox{whenever }\;j=0,\ldots,j_\tau(k)-1\;\mbox{ and }\;j=j^\tau(k)+1,\ldots,k.
\end{array}\right.
\end{equation}
\begin{equation}\label{e:75}
-p^{xk}_{k}\in\lambda^{k}\partial\varphi(\ox^{k}_{k})+\sum^m_{i=1}\alpha^k_ix^*_i\;p^{uk}_{k}=\sum^m_{i=1}\alpha^k_ix^*_i-2\xi^k_k\ou^k_k,\;p^{ak}_{k}=0,
\end{equation}
\begin{equation}\label{e:76}
p^{uk}_{j+1}=\lambda^{k}(v^{uk}_j+h^{-1}_k\theta^{uk}_j),\;p^{ak}_{j+1}=\lambda^{k}(v^{ak}_j+h^{-1}_k\theta^{ak}_j),\quad j=0,\ldots,k-1,
\end{equation}
\begin{equation}\label{e:77}
\begin{array}{ll}
&\bigg(\dfrac{p^{xk}_{j+1}-p^{xk}_{j}}{h_k}-\lambda^kw^{xk}_j-\chi^k_j,\dfrac{p^{uk}_{j+1}-p^{uk}_{j}}{h_k}-\lambda^kw^{uk}_j,\dfrac{p^{ak}_{j+1}-p^{ak}_{j}}{h_k}-\lambda^kw^{ak}_j,\\
&p^{xk}_{j+1}-\lambda^k(v^{xk}_j+h^{-1}_k\theta^{xk}_j)\bigg)\in\left(0,\dfrac{2}{h_k}\xi^k_j\ou^k_j,0,0\right)+N\bigg(\bigg(\ox^{k}_{j},\ou^{k}_{j},\oa^{k}_{j},
\dfrac{\ox^{k}_{j+1}-\ox^{k}_{j}}{-h_k}\bigg);\gph F\bigg)
\end{array}
\end{equation}
for $j=1,\ldots,k-1$ and with the subgradient vectors
\begin{equation}\label{e:78}
\big(w^{xk}_{j},w^{uk}_{j},w^{ak}_{j},v^{xk}_{j},v^{uk}_{j},v^{ak}_{j}\big)\in\partial\ell\bigg(\oz^{k}_{j},\dfrac{\oz^{k}_{j+1}-\oz^{k}_{j}}{h_k}\bigg),\quad j=0,\ldots,k-1,
\end{equation}
where the sequence of $\ve_k\dn 0$ as $k\to\infty$ is taken from \eqref{e:39}, where
\begin{equation}\label{chi}
\chi^k_j:=\left\{\begin{array}{ll}
\disp\frac{x^k_1-x^k_0}{h_k}&\mbox{if }\;j=j=0,\ldots,j_\tau(k)-1\;\mbox{ or }\;j=j^\tau(k)+1,\ldots,k,\\
0&\mbox{otherwise},
\end{array}\right.
\end{equation}
and where the vector triples $(\theta^{xk}_j,\theta^{uk}_j,\theta^{ak}_j)$ for each $j=0,\ldots,k-1$ are defined by
\begin{equation}\label{e:79}
(\theta^{xk}_j,\theta^{uk}_j,\theta^{ak}_j):= 2\int^{t_{j+1}}_{t_j}\bigg(\dfrac{\ox^{k}_{j+1}-\ox^{k}_{j}}{h_k}-\dot{\ox}(t),\dfrac{\ou^{k}_{j+1}-\ou^{k}_{j}}{h_k}-\dot{\ou}(t),\dfrac{\oa^{k}_{j+1}-\oa^{k}_{j}}
{h_k}-\dot{\oa}(t)\bigg)dt.
\end{equation}
\end{theorem}\vspace*{-0.05in}
{\bf Proof.}  Let $y:=(x_0,\ldots,x_k,u_0,\ldots,u_k,a_0,\ldots,a_k,X_0,\ldots,X_{k-1},U_0,\ldots,U_{k-1},A_0,\ldots,A_{k-1})$, where the the starting point $x_0$ is fixed but the other variables depend on $k$ while we omit the upper index ``$k$" for simplicity. Given $\epsilon >0$ in $(P^\tau_k)$, define the problem of mathematical programming $(MP)$ by:
\begin{eqnarray*}
\begin{aligned}
\mbox{minimize }\;&\varphi_{0}[y]:=\varphi(x_{k})+h_{k}\sum^{k-1}_{j=0}\ell(x_j,u_j,a_j,X_j,U_j,A_j)+\sum^{k-1}_{j=0}\int^{t_{j+1}}_{t_j}\left\|(X_j,U_j,A_j)-\dot{\oz}(t)
\right\|^{2}dt\\+&
\left\|\frac{x^k_1-x^k_0}{h_k}-\dot\ox(0)\right\|^2+\dist^2\left(\left\|\dfrac{u^k_1-u^k_0}{h_k}\right\|;(-\infty,\Tilde{\mu}]\right)+\dist^2\left(\sum^{k-2}_{j=0}\left\|U_{j+1}-U_j
\right\|;\big(-\infty,\Tilde{\mu}\big]\right)
\end{aligned}
\end{eqnarray*}
subject to finitely many equality, inequality, and geometric constraints
\begin{eqnarray*}
\begin{aligned}
b^x_j(y):=&x_{j+1}-x_j-h_kX_j=0\;\mbox{ for }\;j=0,\ldots,k-1,\\
b^u_j(y):=&u_{j+1}-u_j-h_kU_j=0\;\mbox{ for }\;j=0,\ldots,k-1,\\
b^a_j(y):=&a_{j+1}-a_j-h_kA_j=0\;\mbox{ for }\;j=0,\ldots,k-1,\\
g_i(y):=&\left\la x^*_i,x_k-u_k\right\ra\le 0\;\mbox{ for }\;i=1,\ldots,m\\
d_j(y):=&\|u_j\|^2-r^2=0\;\mbox{ for }\;j=j_\tau(k),\ldots,j^\tau(k),\\
y\in\Omega_j:=&\big\{y\big|\;r-\tau-\varepsilon_k\le\|u_j\|\le r+\tau+\varepsilon_k\big\}\;\mbox{ for }\;j=0,\ldots,j_\tau(k)-1\;\mbox{ and }\;j=j^\tau(k)+1,\ldots,k,\\
\phi_{j}(y):=&\left\|(x_j,u_j,a_j)-\oz(t_j)\right\|-\epsilon/2\le0\;\mbox{ for }\;j=0,\ldots,k,\\
\phi_{k+1}(y):=&\sum^{k-1}_{j=0}\int^{t_{j+1}}_{t_j}\bigg(\left\|(X_j,U_j,A_j)-\dot{\oz}(t)\right\|^{2}\bigg)dt-\frac{\epsilon}{2}\le0,\\
\phi_{k+2}(y):=&\sum^{k-2}_{j=0}\left\|U_{j+1}-U_j\right\|\le\Tilde{\mu}+1,\quad\phi_{k+3}(y):=\|u_1-u_0\|\le(\Tilde{\mu}+1)(t^k_1-t^k_0),\\
y\in\Xi_{j}:=&\big\{y|-X_{j}\in F(x_{j},u_{j},a_{j})\big\}\;\mbox{ for }\;j=0,\ldots,k-1,\\
y\in\Xi_{k}:=&\big\{y|x_{0}\;\mbox{ is fixed},\;(u_0,a_0)=\big(\ou(0),\oa(0)\big)\big\},
\end{aligned}
\end{eqnarray*}
where the number $\Tilde\mu$ is is calculated in \eqref{e:31}. It follows directly from the construction above that problem $(MP)$ is equivalent to $(P^\tau_k)$ for any fixed $k\in\N$ and $\tau\in[0,\Bar\tau]$.\vspace*{-0.05in}

Necessary optimality conditions for problem $(MP)$ in terms of the generalized differential tools of Section~6 can be deduced from \cite[Theorem~5.24]{m-book2}. We specify them for the optimal solution $\oy=(\oz,\bar{Z})$ to $(MP)$, where $\oz:=(\ox_0,\ldots,\ox_k,\ou_0,\ldots,\ou_k,\oa_0,\ldots,\oa_k)$ is generated by the optimal solution $\oz^{k}$ to $(P^\tau_k)$ while $\bar{Z}:=(\bar{X}_0,\ldots,\bar{X}_{k-1},\bar{U}_0,\ldots,\bar{U}_{k-1},\bar{A}_0,\ldots,\bar{A}_{k-1})$ signifies the discrete ``velocity" determined by the constraints $b_{j}(\oy)=0$. It follows from Theorem~\ref{Th:5} that all the inequality constraints in $(MP)$ relating to functions $\phi_j$ as $j=0,\ldots,k+2$ are inactive for large $k$, and so the corresponding multipliers do not appear in the optimality conditions. Thus we find $\lambda\ge 0,\;\alpha=(\al_1,\ldots,\al_m)\in\R^m_+,\;\xi=(\xi_0,\ldots,\xi_k)\in\R^{k+1},\;p_j=(p^x_j,p^u_j,p^a_j)\in\R^{2n+d}$ as $j=0,\ldots,k$, and
$$
y^*_j=(x^*_{0j},\ldots,x^*_{kj},u^*_{0j},\ldots,u^*_{kj},a^*_{0j},\ldots,a^*_{kj},X^*_{0j},\ldots,X^*_{(k-1)j},U^*_{0j},
\ldots,U^*_{(k-1)j},A^*_{0j},\ldots,A^*_{(k-1)j})
$$
for $j=0,\ldots,k$, which are not zero simultaneously while satisfying \eqref{e:74} and the conditions
\begin{equation}\label{e:80}
y^*_j\in\left\{\begin{array}{ll}
N(\oy;\Xi_j)+N(\oy;\Omega_j)&\mbox{if }\;j\in\big\{0,\ldots,j_\tau(k)-1\big\}\cup\big\{j^\tau(k)+1,\ldots,k\big\},\\
N(\oy;\Xi_j)&\mbox{if }\; j\in\big\{j_\tau(k),\ldots,j^\tau(k)\big\},
\end{array}\right.
\end{equation}
\begin{equation}\label{e:81}
-y^*_0-\ldots-y^*_k\in\lambda\partial\varphi_{0}(\oy)+\sum^m_{i=1}\alpha_i\nabla g_i(\oy)+\sum^{j^\tau(k)}_{j=j_\tau(k)}\xi_j\nabla d_j(\oy)+\sum^{k-1}_{j=0}
\nabla b_j(\oy)^*p_{j+1},
\end{equation}
\begin{equation}\label{e:82}
\alpha_i g_i(\oy)=0\;\mbox{ for }\;i=1,\ldots,m.
\end{equation}
Note that the first line in \eqref{e:80} comes from applying the normal cone intersection formula from \cite[Corollary~3.5]{m-book1} to $\oy\in\Omega_j\cap\Xi_j$ for $j\in \{0,\ldots,j_\tau(k)-1\}\cup\{j^\tau(k)+1,\ldots,k\}$, where the qualification condition imposed therein can be easily verified. It follows from the structure of the sets $\O_j$ and $\Xi_j$ that \eqref{e:74} holds while the inclusions in \eqref{e:80} are equivalent to
\begin{equation}\label{e:83}
\left\{\begin{array}{ll}
\left(x^*_{jj},u^*_{jj},a^*_{jj},-X^{*}_{jj}\right)\in N\left(\left(\ox^{k}_{j},\ou^{k}_{j},\oa^{k}_{j},\dfrac{\ox^{k}_{j+1}-\ox^{k}_{j}}{-h_k}\right);\gph F\right),\;j=j_\tau(k),\ldots,j^\tau(k),\\
\left(x^*_{jj},u^*_{jj}-\psi^u_j,a^*_{jj},-X^{*}_{jj}\right)\in N\left(\left(\ox^{k}_{j},\ou^{k}_{j},\oa^{k}_{j},\dfrac{\ox^{k}_{j+1}-\ox^{k}_{j}}{-h_k}\right);\gph F\right),\;j\not\in\{j_\tau(k),\ldots,j^\tau(k)\}
\end{array}\right.
\end{equation}
with every other components of $y^*_j$ equal to zero, where $\psi^u_j\in N(\ou_j,h^{-1}([r-\tau-\ve_k,r+\tau+\ve_k]))$ and $h(u):=\|u\|^2$. It then follows from \cite[Theorem~1.17]{m-book1} that
$$
N(\ou_j,h^{-1}([r-\tau-\ve_k,r+\tau+\ve_k]))=2\ou_jN(\|\ou_j\|;[r-\tau-\ve_k,r+\tau+\ve_k]).
$$
Hence $\psi^u_j=2\xi_j\ou_j$, where $\xi_j$ satisfies \eqref{e:74} and $j\not\in\{j_\tau(k),\ldots,j^\tau(k)\}$. Similarly we conclude that the triple $(x^*_k(0),u^*_k(0),a^*_k(0))$ determined by the normal cone to  $\Xi_k$ is the only potential nonzero component of $y^*_k$. This shows that
\begin{equation}\label{e:84}
\begin{aligned}
-y^*_0-y^*_1-\ldots-y^*_k=\big(&-x^*_{0k}-x^*_{00},-x^*_{11},\ldots,-x^*_{k-1,k-1},0,-u^*_{0k}-u^*_{00},\ldots,-u^*_{k-1,k-1},0,\\
&-a^*_{0k}-a^*_{00},-a^*_{11},\ldots,-a^*_{k-1,k-1},0,-X^*_{00},\ldots,-X^*_{k-1,k-1},0,\ldots,0\big).
\end{aligned}
\end{equation}
Next we calculate the three sums on the right-hand side of \eqref{e:81}. It is easy to see that
\begin{eqnarray*}
\begin{aligned}
\left(\sum^m_{i=1}\alpha_i\nabla g_i(\oy)\right)_{(x_k,u_k,a_k)}=&\left(\sum^m_{i=1}\alpha_ix^*_i,-\sum^m_{i=1}\alpha_ix^*_i,0\right),\\
\left(\sum^{j^\tau(k)}_{j=j_\tau(k)}\xi_j\nabla d_j(\oy)\right)_{u_j}=&2\xi_j\ou_j\;\mbox{ for }\;j=0,\ldots,k\;\mbox{ with }\;\xi_j=0\;\mbox{ if }\ j\not\in\{j_\tau(k),\ldots,j^\tau(k)\};\\
\left(\sum^{k-1}_{j=0}(\nabla f_j(\oy))^*p_{j+1}\right)_{(x_j,u_j,a_j)}=&
\left\{\begin{array}{lcl}-p_1\quad\mbox{if }\;j=0,\\
p_j-p_{j+1}\quad\mbox{if }\;j=1,\ldots,k-1,\\
p_k\quad\mbox{if }\;j=k,
\end{array}\right.\\
\left(\sum^{k-1}_{j=0}(\nabla f_j(\oy))^*p_{j+1}\right)_{(X,U,A)}=&-h_k p=\big(-h_kp^x_1,\ldots,-h_kp^x_k,-h_kp^u_1,\ldots,-h_kp^u_k,-h_kp^a_1,\ldots,-h_kp^a_k\big).
\end{aligned}
\end{eqnarray*}
Furthermore, the subdifferential sum rule from \cite[Theorem~2.13]{m-book1} gives us the inclusion
$$
\partial\varphi_0(\oy)\subset\partial\varphi(\ox_k)+2\Big(\frac{\ox^k_1-\ox^k_0}{h_k}-\dot\ox(0)\Big)+h_k\sum^{k-1}_{j=0}\partial\ell(\ox_j,\ou_j,\oa_j,\bar{X}_j,
\bar{U}_j,\bar{A}_j)+\sum^{k-1}_{j=0}\nabla
\rho_j(\oy)+\partial\sigma(\oy),
$$
where the functions $\rho(\cdot)$ and $\sigma(\cdot)$ are defined by
$$
\rho_j(y):=\int^{j+1}_j\left\|(X_j,U_j,A_j)-\dot{\oz}(t)\right\|^2dt\quad\mbox{and}
$$
$$
\sigma(y):=\dist^2\left(\left\|\dfrac{u^k_1-u^k_0}{h_k}\right\|;\big(-\infty,\Tilde{\mu}\big]\right)+\dist^2\left(\sum^{k-2}_{j=0}\left\|U_{j+1}-U_j
\right\|;\big(-\infty,\Tilde{\mu}\big]\right).
$$
Note that the function $\psi(x):=\dist^2(x;(-\infty,\Tilde{\mu}])$ is obviously differentiable with $\nabla\phi(x)=0$ for all $x\le\Tilde{\mu}$. Combining this with second estimate in \eqref{e:50} yields $\partial\sigma(\oy)=\{0\}$. Observe also that the nonzero part of $\nabla\rho_j(\oy)$ is calculated by $\nabla_{X_j,U_j,A_j}\rho(\oy)=(\theta^x_j,\theta^u_j,\theta^a_j)$, where the latter triple is defined in \eqref{e:79}. Hence the set $\lambda\partial\varphi_0(\oy)$ in \eqref{e:81} is represented as the collection of
$$
\begin{array}{ll}\lambda(h_kw^x_0,h_kw^k_1+\chi^k_1,\ldots,h_kw^x_{k-1},\vartheta^k,h_kw^u_0,\ldots,h_kw^u_{k-1},0,h_kw^a_0,\ldots,h_kw^a_{k-1},0,\theta^x_0+h_kv^x_0,\ldots,\\
\theta^x_{k-1}+h_kv^x_{k-1},\theta^u_0+h_kv^u_0,\ldots,\theta^u_{k-1}+h_kv^u_{k-1},\theta^a_0+h_kv^a_0,\ldots,\theta^a_{k-1}+h_kv^a_{k-1})
\end{array}
$$
where $\vartheta^k\in\partial\varphi(\ox_k)$, $\chi^k_1$ is defined in \eqref{chi}, and the components of $(w^x,w^u,w^a,v^x,v^u,v^a)$ satisfy \eqref{e:78}. This together with \eqref{e:84} and the above gradient formulas shows that \eqref{e:81} amounts to the relationships
\begin{equation}\label{e:85}
-x^*_{0k}-x^*_{00}=\lambda h_kw^x_0-p^x_1,
\end{equation}
\begin{equation}\label{e:86}
-x^*_{jj}=\lambda h_kw^x_j+\chi^k_j+p^x_j-p^x_{j+1}\;\mbox{ for }\;j=1,\ldots,k-1\;\mbox{ with }\chi_j=0\;\mbox{ if }\;j\ne 1,
\end{equation}
\begin{equation}\label{e:87}
0=\lambda\vartheta^k+\sum^m_{i=1}\alpha_i x^*_i+p^x_k,
\end{equation}
\begin{equation}\label{e:88}
-u^*_{0k}-u^*_{00}=\lambda h_kw^u_0+2\xi_0\ou_0-p^u_1,
\end{equation}
\begin{equation}\label{e:89}
-u^*_{jj}=\lambda h_kw^u_j+2\xi_j\ou_j+p^u_j-p^u_{j+1},
\end{equation}
\begin{equation}\label{e:90}
0=-\sum^k_{i=1}\alpha_i x^*_i+p^u_k+2\xi_k\ou_k,
\end{equation}
\begin{equation}\label{e:91}
-a^*_{0k}-a^*_{00}=\lambda h_kw^a_0-p^a_1,
\end{equation}
\begin{equation}\label{e:92}
-a^*_{jj}=\lambda h_kw^a_j+p^a_j-p^a_{j+1}\;\mbox{ for }\;j=1,\ldots,k-1,
\end{equation}
\begin{equation}\label{e:93}
0=p^a_k,
\end{equation}
\begin{equation}\label{e:94}
-X^*_{jj}=\lambda(h_kv^x_j+\theta^x_j)-h_kp^x_{j+1}\;\mbox{ for }\;j=0,\ldots,k-1,
\end{equation}
\begin{equation}\label{e:95}
0=\lambda(h_kv^u_j+\theta^u_j)-h_kp^u_{j+1}\;\mbox{ for }\;j=0,\ldots,k-1,
\end{equation}
\begin{equation}\label{e:96}
0=\lambda(h_kv^a_j+\theta^a_j)-h_kp^a_{j+1}\;\mbox{ for }\;j=0,\ldots,k-1.
\end{equation}\vspace*{-0.2in}

Now we are ready to justify all the conditions claimed in the theorem. Observe first that \eqref{e:82} clearly yields \eqref{e:73}. Next we extend the vector $p$ by a zero component by putting $p_0:=(x^*_{0k},u^*_{0k},a^*_{0k})$. Then the conditions in \eqref{e:75} follow from \eqref{e:87}, \eqref{e:90}, and \eqref{e:93}. Furthermore, the conditions in \eqref{e:76} follow from \eqref{e:95} and \eqref{e:96}. Using the relationships
$$
\dfrac{p^{x}_{j+1}-p^{x}_{j}}{h_k}-\lambda w^{x}_{j}-\chi^k_j=\dfrac{x^*_{jj}}{h_k},\quad\dfrac{p^{u}_{j+1}-p^{u}_{j}}{h_k}-\lambda w^{u}_{j}=\dfrac{u^{*}_{jj}}{h_k}+2\xi_j\ou_j,
$$
$$
\dfrac{p^{a}_{j+1}-p^{a}_{j}}{h_k}-\lambda w^{a}_{j}=\dfrac{a^{*}_{jj}}{h_k},\quad p^{x}_{j+1}-\dfrac{1}{h_k}\lambda\theta_{xj}-\lambda v^{xk}_{j}=\dfrac{X^{*}_{jj}}{h_k},
$$
which hold due to \eqref{e:86}, \eqref{e:89}, \eqref{e:92}, and \eqref{e:94}, and then substituting them into the left-hand side of \eqref{e:83}, we arrive at \eqref{e:77} for all $j=0,\ldots,k-1$. To verify the nontriviality condition \eqref{e:72}, suppose by contradiction that $\lm=0,\;\alpha=0,\;\xi=0,\; p^{uk}_0=0,\;p^{ak}_0=0$, and $p^x_j=0$ for $j=0,\ldots,k-1$. Then \eqref{e:87} yields that $p^x_k=0$ and thus $p^x_j=0$ for all $j=0,\ldots,k$. Observe further that $x^*_{0k}=p^x_0=0$, and so the conditions in \eqref{e:85}, \eqref{e:86}, and \eqref{e:94} imply that $x^*_{jj}=0$ and $X^*_{jj}=0$ for $j=0,\ldots,k-1$. The validity of \eqref{e:95} and \eqref{e:96} ensures that $p^u_j=0$ and $p^a_j=0$ for $j=1,\ldots,k$, which in turn show by \eqref{e:88}, \eqref{e:89}, \eqref{e:91}, and \eqref{e:92} that $u^*_{jj}=0$ and $a^*_{jj}=0$ for $j=0,\ldots,k-1$. As mentioned above, all the components of $y^*_j$ different from $(x^*_{jj},u^*_{jj},a^*_{jj},X^*_{jj})$ are zero for $j=0,\ldots,k-1$. Hence we have $y^*_j=0$ for $j=0,\ldots,k-1$ and similarly $y^*_k=0$ since the only potential nonzero component of this vector is $x^*_{0k}=p^x_0=0$. We get therefore that $y^*_j=0$ for all $j=0,\ldots,k$, which violates the nontriviality condition for $(MP)$ and thus completes the proof of the theorem. $\h$\vspace*{-0.05in}

The final result of this section employs the effective coderivative calculations for the sweeping control mapping taken from Theorem~\ref{Th:9} that allows us to obtain necessary optimality conditions in $(P^\tau_k)$ expressed entirely via the given problem data and the minimizer under consideration under the additional assumption on the smoothness of $f$, which is imposed for simplicity. Furthermore, we derive an enhanced nontriviality relation in the case of linear independence of the generating vectors $x^*_i$ for the underlying convex polyhedron $C$ from \eqref{e:6}.\vspace*{-0.1in}

\begin{theorem}{\bf (optimality conditions for discretized sweeping inclusions via their initial data).}\label{Th:11} Let $\oz^k=(\ox^k,\ou^k,\oa^k)$ be an optimal solution to problem $(P^\tau_k)$ in the general framework of Theorem~{\rm\ref{Th:10}} with $F$ given by \eqref{e:20} via the active constraint indices $I(\cdot)$ in \eqref{e:19} and locally smooth perturbation function $f$, let the active index subsets $I_0(\cdot)$ and $I_>(\cdot$) be taken from \eqref{e:64}, and let the triples $(\theta^{xk}_j,\theta^{uk}_j,\theta^{ak}_j)$ be defined in \eqref{e:79}. Then there exist dual elements $(\lm^k,\xi^k,p^k)$ as in Theorem~{\rm\ref{Th:10}} together with vectors $\eta^k_j\in\R^m_+$ as $j=0,\ldots,k$ and $\gamma^k_{j}\in\R^m$ as $j=0,\ldots,k-1$ satisfying \eqref{e:74}, the {\sc nontriviality condition}
\begin{equation}\label{e:72a}
\lambda^{k}+\left\|\eta^{k}_k\right\|+\left\|\xi^{k}\right\|+\sum^{k-1}_{j=0}\left\|p^{xk}_{j}\right\|+\left\|p^{uk}_0\right\|+\left\|p^{ak}_0\right\|\not=0,
\end{equation}
the {\sc primal-dual dynamic equations} for all $j=0,\ldots,k-1$ with $\chi^k_j$ defined in \eqref{chi}:
\begin{equation}\label{e:98}
\dfrac{\ox^k_{j+1}-\ox^k_j}{-h_k}-f(\ox^k_j,\oa^k_j)=\sum_{i\in I(\ox^k_j-\ou^k_j)}\eta^k_{ji}x^*_i,
\end{equation}
\begin{equation}\label{e:101}
\begin{array}{ll}
\dfrac{p^{xk}_{j+1}-p^{xk}_j}{h_k}-\lm^kw^{xk}_j-\chi_j^k=\nabla_xf(\ox^k_j,\oa^k_j)^*\left(\lm^k(v^{xk}_j+h^{-1}_k\theta^{xk}_j)-p^{xk}_{j+1}\right)\\
+\disp\sum_{i\in I_0\left(-p^{xk}_{j+1}+\lm^k(h^{-1}_k\theta^k_{xj}+v^{xk}_j)\right)\cup I_>\left(-p^{xk}_{j+1}+\lm^k(h^{-1}_k\theta^k_{xj}+v^{xk}_j)\right)}\gamma^k_{ji}x^*_i,
\end{array}
\end{equation}
\begin{equation}\label{e:102}
\begin{array}{ll}
\dfrac{p^{uk}_{j+1}-p^{uk}_j}{h_k}-\lm^kw^{uk}_j-\dfrac{2}{h_k}\xi^k_j\ou^k_j=-\disp\sum_{i\in I_0\left(-p^{xk}_{j+1}+\lm^k(h^{-1}_k\theta^k_{xj}+v^{xk}_j)\right)\cup I_>\left(-p^{xk}_{j+1}+\lm^k(h^{-1}_k\theta^k_{xj}+v^{xk}_j)\right)}\gamma^k_{ji}x^*_i,
\end{array}
\end{equation}
\begin{equation}\label{e:103}
\dfrac{p^{ak}_{j+1}-p^{ak}_j}{h_k}-\lm^kw^{ak}_j=\nabla_af(\ox^k_j,\oa^k_j)^*\left(\lm^k(v^{xk}_j+h^{-1}_k\theta^{xk}_j)-p^{xk}_{j+1}\right)
\end{equation}
with $(w^{xk}_j,w^{uk}_j,w^{ak}_j,v^{xk}_j,v^{uk}_j,v^{ak}_j)$ taken from \eqref{e:78}, and the right endpoint {\sc transversality conditions}
\begin{equation}\label{e:106a}
-p^{xk}_k\in\lm^k\partial\varphi(\ox^k_k)+\sum^m_{i=1}\eta^k_{ki}x^*_i,\quad p^{uk}_k=\sum^m_{i=1}\eta^k_{ki}x^*_i-2\xi^k_k\ou^k_k,\quad p^{ak}_k=0
\end{equation}
such that the following implications hold:
\begin{equation}\label{e:97}
\big[\left\la x^*_i,\ox^k_j-\ou^k_j\right\ra<0\big]\Longrightarrow\eta^k_{ji}=0,
\end{equation}
\begin{equation}\label{e:100}
\left\{\begin{array}{lcl}
i\in I_0\left(-p^{xk}_{j+1}+\lm^k(h^{-1}_k\theta^k_{xj}+v^{xk}_j)\right)\Longrightarrow\gamma^k_{ji}\in\R,\\
i\in I_>\left(-p^{xk}_{j+1}+\lm^k(h^{-1}_k\theta^k_{xj}+v^{xk}_j)\right)\Longrightarrow\gamma^k_{ji}\ge 0,\\
\big[i\not\in I_0\left(-p^{xk}_{j+1}+\lm^k(h^{-1}_k\theta^k_{xj}+v^{xk}_j)\right)\cup I_>\left(-p^{xk}_{j+1}+\lm^k(h^{-1}_k\theta^k_{xj}+v^{xk}_j)\right)\big] \Longrightarrow\gamma^k_{ji}=0
\end{array}\right.
\end{equation}
for $j=0,\ldots,k-1$ and $i=1,\ldots,m$. Furthermore, we have the constraint conditions \eqref{e:74} together with
\begin{equation}\label{e:104}
\big[\la x^*_i,\ox^k_j-\ou^k_j\ra<0\big]\Longrightarrow\gamma^k_{ji}\;\mbox{ for }\;j=0,\ldots,k-1\;\mbox{ and }\;i=1,\ldots,m,
\end{equation}
\begin{equation}\label{e:106}
\big[\la x^*_i,\ox^k_k-\ou^k_k\ra<0\big]\Longrightarrow\eta^k_{ki}=0\;\mbox{ for }\;i=1,\ldots,m.
\end{equation}
Finally, the linear independence of the vectors $\{x^*_i|\;i\in I(\ox^k_j-\ou^k_j)\}$ in \eqref{e:6} ensures the implication
\begin{equation}\label{e:99}
\eta^k_{ji}>0\Longrightarrow\big[\la x^*_i,-p^{xk}_{j+1}+\lm^k\left(h^{-1}_k\theta^k_{xj}+v^{xk}_j\right)\ra=0\big]
\end{equation}
and the validity of the {\sc enhanced nontriviality condition}
\begin{equation}\label{e:108}
\lm^k+\left\|\xi^k\right\|+\left\|p^{uk}_0\right\|+\left\|p^{ak}_0\right\|\not=0.
\end{equation}
\end{theorem}\vspace*{-0.05in}
{\bf Proof.} The coderivative definition \eqref{cod} allows us to equivalently rewrite the discrete Euler-Lagrange inclusion \eqref{e:77} in the coderivative form
\begin{equation}\label{e:109}
\begin{array}{l}
\bigg(\dfrac{p^{xk}_{j+1}-p^{xk}_j}{h_k}-\lm^kw^{xk}_j,\dfrac{p^{uk}_{j+1}-p^{uk}_j}{h_k}-\lm^kw^{uk}_j-\dfrac{2}{h_k}\xi^k_j\ou^k_j,\dfrac{p^{ak}_{j+1}-p^{ak}_j}
{h_k}-\lm^kw^{ak}_j\bigg)\\
\in D^*F\bigg(\ox^k_j,\ou^k_j,\oa^k_j,\dfrac{\ox^k_{j+1}-\ox^k_j}{-h_k}\bigg)\left(\lm^k(h^{-1}_k\theta^k_{xj}+v^{xk}_j)-p^{xk}_{j+1}\right),\quad j=0,\ldots,k-1.
\end{array}
\end{equation}
Taking into account that $\dfrac{\ox^k_{j+1}-\ox^k_j}{-h_k}-f\left(\ox^k_j,\oa^k_j\right)\in G\left(\ox^k_j-\ou^k_j\right)$ for $j=0,\ldots,k-1$ with $G(x)=N(x;C)$, we find by \eqref{e:109} vectors $\eta^k_j\in\R^m_+$ for $j=0,\ldots,k-1$ such that conditions \eqref{e:98} and \eqref{e:97} hold. Using now the coderivative inclusion \eqref{e:70} from Theorem~\ref{Th:9} with $x:=\ox^k_j,\;u:=\ou^k_j,\;a:=\oa^k_j,\;w:=\dfrac{\ox^k_{j+1}-\ox^k_j}{-h_k}$, and $y:=\lm^k(h^{-1}_k\theta^k_{xj}+v^{xk}_j)-p^{xk}_{j+1}$ shows $\gamma^k_j\in\R^m$ and the relationships
\begin{equation*}
\begin{aligned}
&\bigg(\dfrac{p^{xk}_{j+1}-p^{xk}_j}{h_k}-\lm^kw^{xk}_j,\dfrac{p^{uk}_{j+1}-p^{uk}_j}{h_k}-\lm^kw^{uk}_j-\dfrac{2}{h_k}\xi^k_j\ou^k_j,\dfrac{p^{ak}_{j+1}-p^{ak}_j}{h_k}-\lm^kw^{ak}_j\bigg)\\
=&\bigg(\nabla_xf(\ox^k_j,\oa^k_j)^*\left(\lm^k(v^{xk}_j+h^{-1}_k\theta^{xk}_j)-p^{xk}_{j+1}\right)\\
+&\sum_{i\in I_0\left(-p^{xk}_{j+1}+\lm^k(h^{-1}_k\theta^k_{xj}+v^{xk}_j)\right)\cup I_>\left(-p^{xk}_{j+1}+\lm^k(h^{-1}_k\theta^k_{xj}+v^{xk}_j)\right)}\gamma^k_{ji}x^*_i,\\
-&\sum_{i\in I_0\left(-p^{xk}_{j+1}+\lm^k(h^{-1}_k\theta^k_{xj}+v^{xk}_j)\right)\cup I_>\left(-p^{xk}_{j+1}+\lm^k(h^{-1}_k\theta^k_{xj}+v^{xk}_j)\right)}\gamma^k_{ji}x^*_i,\nabla_af(\ox^k_j,\oa^k_j)^*\left(\lm^k(v^{xk}_j+h^{-1}_k\theta^{xk}_j)
-p^{xk}_{j+1}\right)\bigg)
\end{aligned}
\end{equation*}
are satisfied for all $j=0,\ldots,k-1$ and thus ensure the validity of all the conditions in \eqref{e:101}, \eqref{e:102}, \eqref{e:103}, \eqref{e:100}, and \eqref{e:104}. Defining now $\eta^k_k:=\alpha_k$ via $\alpha_k$ from the statement of Theorem~\ref{e:72} yields $\eta^k_j\in\R^m_+$ for $j=0,\ldots,k$ and allows us to deduce the nontriviality condition \eqref{e:72a} from that in \eqref{e:72} and the transversality conditions in \eqref{e:106a} from those in \eqref{e:75} and \eqref{e:76}. Implication \eqref{e:106} is a direct consequence of \eqref{e:73} and the definition of $\eta^k_k$.\vspace*{-0.05in}

Assume finally that the generating vectors  $\{x^*_i|\;i\in I(\ox^k_j-\ou^k_j)\}$ of the convex polyhedron $C$ are linearly independent. It is not hard to observe that the inclusion
$$
\lm^k(h^{-1}_k\theta^k_{xj}+v^{xk}_j)-p^{xk}_{j+1}\in\dom D^*G\left(\ox^k_j-\ou^k_j,\dfrac{\ox^k_{j+1}-\ox^k_j}{-h_k}-f(\ox^k_j,\oa^k_j)\right),
$$
which follows from \eqref{e:109}, in this case yields \eqref{e:99}. It remains to verify the enhanced nontriviality condition \eqref{e:108} under the imposed linear independence. Suppose on the contrary that $\lm^k=0,\;\xi^k=0,\;p^{uk}_0=0$, and $p^{ak}_0=0$. Then $p^{uk}_j=0,\;p^{ak}_j=0$ for $j=0,\ldots,k$ by \eqref{e:76}. Furthermore, it follows from the second condition in \eqref{e:106a} with $p^{uk}_k=0$ that $\sum^m_{i=1}\eta^k_{ki}x^*_i=0$. This implies by
definition \eqref{e:19} of the active constraint indices and the imposed linear independence of $x^*_i$ over this index set that $\eta^k_k=0$, and so $p^{xk}_k=0$ by the first condition in \eqref{e:106a}. On the other hand, we get from \eqref{e:102} that
$$
\sum_{i\in I_0\left(-p^{xk}_{j+1}+\lm^k(h^{-1}_k\theta^k_{xj}+v^{xk}_j)\right)\cup I_>\left(-p^{xk}_{j+1}+\lm^k(h^{-1}_k\theta^k_{xj}+v^{xk}_j)\right)}\gamma^k_{ji}x^*_i=0.
$$
Combining this with \eqref{e:101} and $p^{xk}_k=0$ shows that $p^{xk}_j=0$ for all $j=0,\ldots,k-1$. This clearly implies that \eqref{e:72a} is violated and hence justifies the validity of \eqref{e:108}. $\h$\vspace*{-0.15in}

\section{Numerical Example}
\setcounter{equation}{0}\vspace*{-0.1in}

In this section we present a numerical example, which illustrates the scheme of applying the obtained necessary optimality conditions to solve sweeping control problems via their discrete approximations.\vspace*{-0.1in}

\begin{example}{\bf(calculating optimal controls).}\label{Ex:1} Consider the original optimal control problem $(P)$ for the sweeping process with $n=m=d=r=T=1$, $x_0=0$, $x^*_1=1$ and following functional data:
\begin{equation}\label{e:111}
\varphi(x):=\dfrac{(x-1)^2}{2},\quad\ell(t,x,u,a,\dot{x},\dot{u},\dot{a}):=\dfrac{1}{2}a^2,\quad f(x,a):=a.
\end{equation}
Since the $u$-component is not involved in the optimization of \eqref{e:111}, we can put $u(t)\equiv 1$ and do not distinguish problems $(P)$ and $(P^\tau)$. Furthermore, we suppose without lost of generality that $a(\cdot)$ is uniformly bounded and observe that all the assumptions of Theorem~\ref{Th:2} and other results of this paper are satisfied for \eqref{e:111}. Thus this problem admits an optimal solution, and we do not need either global or local relaxation therein. The sweeping dynamics in $(P)$ can be written as
\begin{equation}\label{e:112}
-\dot{x}(t)\in N\big(x(t);C(t)\big)+f\big(a(t)\big)\,\mbox{ a.e. }\;t\in[0,1],
\end{equation}\vspace*{-0.1in}
with
$$
C(t):=\big\{x\in\R\big|\;x\le 0\big\}+1 \ \mbox{for all}\ t\in[0,1].
$$
To start applying necessary optimality conditions of Theorem~\ref{Th:11} for the discrete approximations of $(P)$, suppose further that $\ox(t)\in{\rm int}(C+\ou(t))$ for any $t\in[0,1)$. This shows that $\la x^*_1,\ox^k_j-\ou^k_j\ra<0$, and thus $\eta^k_j=0$ for all $j=0,\ldots,k$. Employing now Theorem~\ref{Th:11} gives us the following relationships:
\begin{enumerate}
\item $w^k_j=(0,0,\oa^k_j)$ and $v^k_j=(0,0,0)$ for $j=0,\ldots,k-1$;
\item $\ox^k_{j+1}-\ox^k_j=-h_k\oa^k_j$ for $j=0,\ldots,k-1$;
\item $\dfrac{p^{xk}_{j+1}-p^{xk}_j}{h_k}=\chi_j^k$ for $j=0,\ldots,k-1$, where $\chi_j^k$ is defined in \eqref{chi};
\item $\dfrac{p^{uk}_{j+1}-p^{uk}_j}{h_k}-\dfrac{2}{h_k}\xi^k_j=0$ for $j=0,\ldots,k-1$;
\item $\dfrac{p^{ak}_{j+1}-p^{ak}_j}{h_k}-\lm^k\oa^k_j=\dfrac{\lm^k\theta^{xk}_j}{h_k}-p^{xk}_{j+1}$ for $j=0,\ldots,k-1$;
\item $\xi^k_j=0$ for $j=0,\ldots,j_\tau(k)-1$ and $j=j^\tau(k)+1,\ldots,k$;
\item $-p^{xk}_k=\lm^k(\ox^k_k-1)$;
\item $p^{uk}_{k}=p^{ak}_{k}=0$;
\item $\lm^k+\|\xi^k\|+|p^{uk}_0|+|p^{ak}_0|\not=0$;
\item $p^{uk}_{j+1}=\dfrac{\lm^k\theta^{uk}_j}{h_k}$ and $p^{ak}_{j+1}=\dfrac{\lm^k\theta^{ak}_j}{h_k}$ for $j=0,\ldots,k-1$;
\item $\theta^{xk}_j=2[(\ox^k_{j+1}-\ox^k_j)-(\ox(t_{j+1})-\ox(t_j))]=2[-h_k\oa^k_j-(\ox(t_{j+1})-\ox(t_j))]$ for $j=0,\ldots,k-1$;
\item $\theta^{uk}_j=2[(\ou^k_{j+1}-\ou^k_j)-(\ou(t_{j+1})-\ou(t_j))]=0$, for $j=0,\ldots,k-1$;
\item $\theta^{ak}_j=2[(\oa^k_{j+1}-\oa^k_j)-(\oa(t_{j+1})-\oa(t_j))]$ for $j=0,\ldots,k-1$.
\end{enumerate}
It follows from the relationships (3) and (7) above that
\begin{equation*}
p^{xk}_j=\left\{
\begin{array}{ll}
-\lm^k(\ox^k_k-1)\;\mbox{ for all}\ j=1,\ldots,k,\\
-\lm^k(\ox^k_k-1)-h_k\chi^k_0\;\mbox {for}\;j=0.
\end{array}\right.
\end{equation*}
Employing now (2) tells us that
$$
\ox^k_{j+1}=-h_k(\oa^k_0+\ldots+\oa^k_j)\;\mbox{ for }\;j=0,\ldots,k-1.
$$
Furthermore, from (5) we have that
$$
\dfrac{p^{ak}_{j+1}-p^{ak}_0}{h_k}=\lm^k(\oa^k_0+\ldots+\oa^k_j)+\lm^k\dfrac{\theta^{xk}_0+\ldots+\theta^{xk}_j}{h_k}-(p^{xk}_1+\ldots+p^{xk}_{j+1}),\quad j=0,\ldots,k-1.
$$
Combining this with (11) brings us to
\begin{equation*}
\begin{aligned}
\dfrac{p^{ak}_{j+1}-p^{ak}_0}{h_k}&=\lm^k(\oa^k_0+\ldots+\oa^k_j)+\lm^k\dfrac{2[-h_k(\oa^k_0+\ldots+\oa^k_j)-(\alpha_1+\ldots+\alpha_{j+1})]}{h_k}+(j+1)\lm^k
(\ox^k_k-1)\\
&= \lm^k(\oa^k_0+\ldots+\oa^k_j)-\dfrac{2\lm^k}{h_k}\big(\alpha_1+\ldots+\alpha_{j+1}\big)+(j+1)\lm^k\big[-h_k(\oa^k_0+\ldots+\oa^k_{k-1})-1\big],
\end{aligned}
\end{equation*}
where $\alpha_j:=\ox(t_j)-\ox(t_{j-1})$ for $j=1,\ldots,k$. Setting $j:=k-1$ and using $p^{ak}_k=0$ give us
\begin{equation}
\label{e:113}
p^{ak}_0=\lm^k+2\lm^k(\alpha_1+\ldots+\alpha_k).
\end{equation}
Thus we arrive at the explicit expression for the adjoint arc as $j=0,\ldots,k-1$:
$$
p^{ak}_{j+1}=\lm^k\big[1-(j+1)h_k\big]+2\lm^k(\alpha_{j+2}+\ldots+\alpha_k)+\lm^kh_k\big[(\oa^k_0+\ldots+\oa^k_j)-(j+1)h_k(\oa^k_0+\ldots+\oa^k_{k-1})\big].
$$
From the latter expression we get the representation
$$
\dfrac{2\lm^k}{h_k}\big(\oa^k_{j+1}-\oa^k_j\big)=\dfrac{2\lm^k}{h_k}\beta_{j+1}+\lm^k\big[1-(j+1)h_k\big]+2\lm^k(\alpha_{j+2}+\ldots+\alpha_k)+\lm^kh_k\big[S_j-(j+1)
h_kS_{k-1}\big],
$$
where $\beta_j:=\oa(t_j)-\oa(t_{j-1})$ for $j=1,\ldots,k$ and $S_j:=\oa^k_0+\ldots+\oa^k_j$ for $j=0,\ldots,k$.\vspace*{-0.05in}

Let us show that $\lm^k\not=0$. Supposing on the contrary that $\lm^k=0$ gives us by (4), (6), (10), and \eqref{e:113} that $\xi^k=0,\;p^{uk}_0=0$, and $p^{ak}_0=0$, which contradicts the nontriviality condition (9). Thus $\lm^k\not=0$ and
\begin{equation*}
\left\{
\begin{array}{ll}
S_{j+1}=\left(1+\dfrac{h^2_k}{2}\right)S_j-h^2_k(j+1)S_{k-1}+\beta_{j+1}+\dfrac{h_k[1-(j+1)h_k]}{2}+h_k(\alpha_{j+2}+\ldots+\alpha_k),\;j=0,\ldots,k-1,\\
S_0=\oa^k_0=\oa(0).
\end{array}\right.
\end{equation*}
Using these equations, we can solve $S_j$ explicitly and so completely calculate $\oa^k_j=S_{j+1}-S_j$ and $\ox^k_j$.
\end{example}\vspace*{-0.2in}

\section{Concluding Remarks}
\setcounter{equation}{0}\vspace*{-0.1in}

The main results of this paper prove the strong $W^{1,2}$-approximation of local optimal solutions to continuous-time sweeping control problems of the hysteresis type by optimal solutions of their well-posed discrete approximations and then establish necessary optimality conditions for discrete optimal solutions entirely via the initial data of the perturbed sweeping process under consideration. The obtained results justify the possibility of solving a new class of highly involved optimal control problems by using finite-dimensional approximations. A challenging issue remains on deriving necessary optimality conditions for local solutions to continuous-time sweeping control problems of this class by passing to the limit from those obtained here for their finite-difference counterparts. Besides their own theoretical interest, explicit necessary optimality conditions for continuous-time sweeping systems may be convenient for calculating optimal solutions. We pursue these goals (in both theory and applications) in our on-going research.

{\bf Acknowledgements.} The authors are indebted to Giovanni Colombo and Nguyen Hoang for helpful discussions on discrete approximations of non-Lipschitzian differential inclusions.
\end{document}